# ASYMPTOTIC MINIMAXITY OF FALSE DISCOVERY RATE THRESHOLDING FOR SPARSE EXPONENTIAL DATA[1]

BY DAVID DONOHO AND JIASHUN JIN

*Stanford University and Purdue University*

We apply FDR thresholding to a non-Gaussian vector whose coordinates $X_i$, $i = 1, \ldots, n$, are independent exponential with individual means $\mu_i$. The vector $\mu = (\mu_i)$ is thought to be sparse, with most coordinates 1 but a small fraction significantly larger than 1; roughly, most coordinates are simply 'noise,' but a small fraction contain 'signal.' We measure risk by per-coordinate mean-squared error in recovering $\log(\mu_i)$, and study minimax estimation over parameter spaces defined by constraints on the per-coordinate $p$-norm of $\log(\mu_i)$: $\frac{1}{n}\sum_{i=1}^{n} \log^p(\mu_i) \leq \eta^p$.

We show for large $n$ and small $\eta$ that FDR thresholding can be nearly Minimax. The FDR control parameter $0 < q < 1$ plays an important role: when $q \leq 1/2$, the FDR estimator is nearly minimax, while choosing a fixed $q > 1/2$ prevents near minimaxity.

These conclusions mirror those found in the Gaussian case in Abramovich et al. [*Ann. Statist.* **34** (2006) 584–653]. The techniques developed here seem applicable to a wide range of other distributional assumptions, other loss measures and non-i.i.d. dependency structures.

**1. Introduction.** Suppose that we have $n$ measurements $X_i$ which are exponentially distributed, with (possibly different) means $\mu_i$:

(1.1) $\qquad X_i \sim \text{Exp}(\mu_i), \qquad \mu_i \geq 1, \ i = 1, \ldots, n.$

The unknown $\mu_i$'s exhibit *sparse heterogeneity*: most take the common value 1, but a small fraction take different values greater than 1.

There are various ways to precisely define sparsity; see [3], for example. In our setting of exponential means, the most intuitive notion of sparsity is

Received May 2004; revised February 2006.
[1]Supported in part by NSF Grants DMS 00-77261 and DMS 95-05151.

simply that there is a *relatively small proportion of $\mu_i$'s which are strictly larger than* 1:

(1.2) $$\frac{\#\{i : \mu_i \neq 1\}}{n} \leq \varepsilon \approx 0.$$

Such situations arise in several areas of application.

- *Multiple lifetime analysis.* Suppose that the $X_i$ represent failure times of many comparable independent systems, where a small fraction of the systems—we do not know which ones—may have significantly higher expected lifetimes than the typical system.
- *Multiple testing.* Suppose that we conduct many independent statistical hypothesis tests, each yielding a *p*-value $p_i$, say, and that the vast majority of those tests correspond to cases where the null distribution is true, while a small fraction correspond to cases where a Lehmann alternative [13] is true. Then $X_i \equiv \log(1/p_i) \sim \text{Exp}(\mu_i)$, where most of the $\mu_i$ are 1, corresponding to true null hypotheses, while a few are greater than 1, corresponding to Lehmann alternatives.
- *Signal analysis.* A common model (e.g., in spread-spectrum communications) for a discrete-time signal $(Y_t)_{t=1}^n$ takes the form $Y_t = \sum_j W_j \times \exp\{\sqrt{-1}\lambda_j t\} + Z_t$, where $Z_t$ is a white Gaussian noise and the $\lambda_j$ index a small number of unknown frequencies with white Gaussian noise coefficients $W_j$. In spectral analysis of such signals, it is common to compute the periodogram $I(\omega) = |n^{-1/2}\sum_t Y_t \exp(\sqrt{-1}\omega t)|^2$ and consider as primary data the periodogram ordinates $X_i \equiv I(\frac{2\pi i}{n})$, $i = 1, \ldots, n/2 - 1$. These can be modeled as independently exponentially distributed with means $\mu_i$, say; here, most of the $\mu_i = 1$, meaning that there is only noise at those frequencies, while some of the $\mu_i > 1$, meaning that there is signal at those frequencies (i.e., certain frequencies $\omega_i = \frac{2\pi i}{n}$ happen to match some $\lambda_j$). In an incoherent or noncooperative setting, we would not know the $\lambda_j$ and, hence, would not know which $\mu_i > 1$.

The simple sparsity model (1.2) is merely a first pass at the problem; in applications, we may also need to consider situations with a large number of means which are close to, but not exactly, 1. A more general assumption (adapted from [3, 7]) is that for some $0 < p < 2$, the log means obey an $\ell^p$ constraint,

$$\frac{1}{n}\left(\sum_{i=1}^n \log^p \mu_i\right) \leq \eta^p, \qquad \eta \text{ small, } 0 < p < 2.$$

Working on the log-scale turns out to be useful because of the 'multiplicative' nature of the exponential data. The parameter $p$ measures the degree of sparsity of $\mu$. As $p \to 0$,

$$\sum_{i=1}^{n} \log^p(\mu_i) \longrightarrow \#\{i : \mu_i \neq 1\}.$$

1.1. *Minimax estimation of sparse exponential means.* We now turn to simultaneous estimation of the means $\mu_i$. Let $\mu = (\mu_1, \mu_2, \ldots, \mu_n)$ and suppose we use the squared $\ell^2$-norm on the log-scale to measure loss,

$$\|\log \hat{\mu} - \log \mu\|_2^2 = \sum_{i=1}^{n} (\log \hat{\mu}_i - \log \mu_i)^2.$$

Motivated by situations of sparsity, we consider restricted parameter spaces, namely $\ell^p$-balls with radius $\eta$,

(1.3) $$M_{n,p}(\eta) = \left\{ \mu : \frac{1}{n} \sum_{i=1}^{n} \log^p(\mu_i) \leq \eta^p \right\}.$$

We quantify performance by means of the expected coordinatewise loss

$$R_n(\hat{\mu}, \mu) = \mathcal{E}\left[ \frac{1}{n} \sum_{i=1}^{n} (\log \hat{\mu}_i - \log \mu_i)^2 \right].$$

We are interested in the minimax risk, the optimal risk which any estimator can guarantee to hold uniformly over the parameter space

(1.4) $$R_n^* = R_n^*(M_{n,p}(\eta)) = \inf_{\hat{\mu}} \sup_{M_{n,p}(\eta)} R_n(\hat{\mu}, \mu).$$

This quantity has been studied before in a related Gaussian noise setting [3], but not, to our knowledge, in an exponential noise setting. Its asymptotic behavior as $\eta \to 0$ is pinned down by the following result:

THEOREM 1.1.
$$\lim_{\eta \to 0} \left[ \frac{\lim_{n \to \infty} R_n^*(M_{n,p}(\eta))}{\eta^p \log^{2-p} \log \frac{1}{\eta}} \right] = 1.$$

A natural approach in this problem is simple thresholding. More precisely, set $\hat{\mu}_t \equiv (\hat{\mu}_{t,i})_{i=1}^n$, where

(1.5) $$\hat{\mu}_{t,i} = \begin{cases} X_i, & X_i \geq t, \\ 1, & \text{otherwise.} \end{cases}$$

For an appropriate choice of threshold $t$ (which depends in principle on $p$ and $\eta$, but not on $n$), this can be asymptotically minimax, as the following result shows:

THEOREM 1.2.
$$\liminf_{\eta \to 0} \left[ \lim_{n \to \infty} \frac{\sup_{M_{n,p}(\eta)} R_n(\hat{\mu}_t, \mu)}{R_n^*(M_{n,p}(\eta))} \right] = 1.$$

Here, by "asymptotically minimax," we mean that the ratio of the worst risk obtained by the estimator to the corresponding minimax risk tends to 1 as $n \to \infty$, followed by $\eta \to 0$.

The minimizing threshold $t_0 = t_0(p, \eta)$ referred to in this theorem behaves as

$$t_0(p, \eta) \sim p \log(1/\eta) + p \log \log(1/\eta) \cdot (1 + o(1)), \qquad \eta \to 0.$$

In order to have asymptotic minimaxity, it is important to adapt the threshold to the sparsity parameters $(p, \eta)$.

1.2. *FDR thresholding.* FDR-controlling methods were first proposed in a multiple hypothesis testing situation in [1, 2]. For the exponential model we are considering, we suppose that there are $n$ independent tests of unrelated hypotheses, $H_{0,i}$ versus $H_{1,i}$, where the test statistics $X_i$ obey the conditions

(1.6) $\qquad$ under $H_{0,i}$: $\qquad X_i \sim \text{Exp}(1)$,

(1.7) $\qquad$ under $H_{1,i}$: $\qquad X_i \sim \text{Exp}(\mu_i)$, $\qquad \mu_i > 1$,

and it is unknown how many of the alternative hypotheses are likely to be true. Select a value $q$, $0 < q < 1$, which Abramovich et al. [1, 2] called the *FDR control parameter*. If we call any case where $H_{0,i}$ is rejected in favor of $H_{1,i}$ a 'discovery,' then a 'false discovery' is a situation where $H_{0,i}$ is falsely rejected. An FDR-controlling procedure controls

$$\mathcal{E}\left[\frac{\#\{\text{False Discoveries}\}}{\#\{\text{Total Discoveries}\}}\right] \le q.$$

Simes' procedure [17] was shown by [4] to be FDR-controlling and it is easy to describe. We begin by sorting all of the observations into descending order,

$$X_{(1)} \ge X_{(2)} \ge \cdots \ge X_{(n)}.$$

Next, compare the sorted values with quantiles of Exp(1); more specifically, if $E(t)$ denotes the standard exponential distribution function and $\bar{E} = 1 - E$ the corresponding survival function, compare $(X_{(1)}, X_{(2)}, \ldots, X_{(n)})$ with $(t_1, t_2, \ldots, t_n)$, where

$$t_k = \bar{E}^{-1}\left(q \cdot \frac{k}{n}\right) = -\log\left(q \cdot \frac{k}{n}\right), \qquad 1 \le k \le n,$$

and let $t_0 = \infty$. Finally, let $k = k_{FDR}$ be the largest index $k \geq 1$ for which $X_{(k)} \geq t_k$, with $k = 0$ if there is no such index. The FDR thresholding estimator $\hat{\mu}_{q,n}^{FDR}$ uses the (data-dependent) threshold $\hat{t}^{FDR} \equiv t_{k_{FDR}}$ and has components $(\hat{\mu}_i)_{i=1}^n$, where

$$\hat{\mu}_i = \begin{cases} X_i, & X_i \geq \hat{t}^{FDR}, \\ 1, & \text{otherwise.} \end{cases} \quad (1.8)$$

In particular, if $k_{FDR} = 0$, then $\hat{\mu}_i = 1$ for all $i$. We think of the observations exceeding $t^{FDR}$ as *discoveries*; the FDR property guarantees relatively few false discoveries.

An attractive property of the procedure is its simplicity and definiteness. Another attractive property is its good performance in an estimation context. Our main result in this paper is the following theorem:

THEOREM 1.3.  1. *When $0 < q \leq \frac{1}{2}$, the FDR estimator $\hat{\mu}_{q,n}^{FDR}$ is asymptotically minimax, that is,*

$$\lim_{\eta \to 0} \left[ \lim_{n \to \infty} \frac{\sup_{\mu \in M_{n,p}(\eta)} R_n(\hat{\mu}_{q,n}^{FDR}, \mu)}{R_n^*(M_{n,p}(\eta))} \right] = 1.$$

2. *When $q > \frac{1}{2}$, the FDR estimator $\hat{\mu}_{q,n}^{FDR}$ is not asymptotically minimax, that is,*

$$\lim_{\eta \to 0} \left[ \lim_{n \to \infty} \frac{\sup_{\mu \in M_{n,p}(\eta)} R_n(\hat{\mu}_{q,n}^{FDR}, \mu)}{R_n^*(M_{n,p}(\eta))} \right] = \frac{q}{1-q} > 1.$$

1.3. *Interpretation.* By controlling the FDR so that *there are at least as many 'true' discoveries above threshold as 'false' ones*, we obtain an estimator that with increasing sparsity $\eta \to 0$, asymptotically attains the minimax risk. This is the case across a wide range of measures of sparsity.

The same general conclusion was found in a model of Gaussian observations due to Abramovich, Benjamini, Donoho and Johnstone [3]. In that setting, the authors supposed that $X_i \sim N(\mu_i, 1)$ and that the $\mu_i$ are mostly close to zero so that $\frac{1}{n}(\sum_{i=1}^n |\mu_i|^p) \leq \eta_n^p$. (Note that the sparsity parameter $\eta$ was replaced by a sequence $\eta_n \to 0$ as $n \to \infty$ in [3].) In that setting, it was shown that FDR thresholding gave asymptotically minimax estimators. Hence, the results in our paper show that FDR thresholding, known previously to be successful in the Gaussian case, is also successful in an interesting non-Gaussian case.

It appears to us that there may be a wide range of non-Gaussian cases wherein the vector of means is sparse and FDR gives nearly-minimax results. Elsewhere, Jin [12] will report results showing that similar conclusions are possible in the case of Poisson data. In that setting, we have, for large $n$, $n$

Poisson observations $N_i \sim \text{Poisson}(\mu_i)$ with the $\mu_i$ mostly 1 and perhaps a small fraction significantly greater than 1. In that setting as well, it seems that FDR thresholding gives near-minimax risk.

In fact, the approach developed here seems applicable to a wide range of non-Gaussian distributions and loss functions. At the same time, it also seems able to cover a wide range of dependence structures.

1.4. *Contents.* The paper is organized as follows. Theorems 1.1 (on minimax risk) and 1.2 (on thresholding risk) are developed and proved in Sections 2 and 3, respectively. These sections also introduce a model in which the parameter $\mu$ is realized by i.i.d. random sampling rather than as a fixed vector; this model is very useful for computations.

Sections 4–7 develop our technical approach to analyzing FDR thresholding. This starts in Section 4 with a definition and analysis of the so-called *FDR functional*, establishing various boundedness and continuity properties. The FDR functional allows us to articulate the idea that in a Bayesian setting where both the mean vector $\mu$ and the subordinate data $X$ are realizations of iid random variables, there is a 'large-sample threshold' which FDR thresholding is consistently 'estimating.' Section 5 discusses the performance of an idealized pseudo-estimator which thresholds at this large-sample threshold even in finite samples; it shows that the idealized 'estimator' achieves risk performance approaching the minimax risk. Section 6 shows that in large samples, the risk of FDR thresholding is well approximated by the risk of idealized FDR thresholding. Section 7 ties together the pieces by showing that the results of Sections 4–6 for the Bayesian model have close parallels in the original frequentist setting of this introduction, implying Theorem 1.3.

Section 8 ends the paper by (i) graphically illustrating two important points about the method and the proof below, (ii) by comparing our results to recent work of Genovese and Wasserman and of Abramovich et al. and (iii) describing generalizations to a variety of non-Gaussian and dependent data structures.

1.5. *Notation.* In this paper, we let $E$ denote the cumulative distribution function (cdf) of $\text{Exp}(1)$, while, to avoid confusion, we use $\mathcal{E}$ for the expectation operator applied to random variables; we also let $\bar{E}$ denote the survival function of $\text{Exp}(1)$ and extend this notation to all cdf's; that is, for any cdf $G$, we let $\bar{G} = 1 - G$ denote the survival function.

We let '#' denote the scale mixture operator, mapping any (marginal) distribution $F$ on $[1, \infty)$ to a corresponding $G = E \# F$ on $[0, \infty)$, according to

$$F \stackrel{E\#}{\longmapsto} G: \qquad G(t) = \int E(t/\mu) \, dF(\mu).$$

Note here that $G$ is the cdf of a scalar random variable $X$, with $\mu$ a random variable $\mu \sim F$ and $X|\mu \sim \text{Exp}(\mu)$. We let $\mathcal{F}$ denote the set of all eligible cdf's,

$$\mathcal{F} = \{F : P_F\{\mu \geq 1\} = 1\},$$

and $\mathcal{F}_p(\eta)$ denote the convex set of $p$th moment-constrained cdf's,

(1.9) $$\mathcal{F}_p(\eta) = \left\{ F \in \mathcal{F} : \int \log^p(\mu)\, dF(\mu) \leq \eta^p \right\}, \qquad 0 < p < 2.$$

We also let $\mathcal{G}$ denote the collection of all scale mixtures of exponentials,

$$\mathcal{G} = \{G : G = E\#F, F \in \mathcal{F}\},$$

and let $\mathcal{G}_p(\eta)$ denote the subclass where the mixing distributions obey the moment condition $\mathcal{E}[\log^p(\mu)] \leq \eta^p$,

(1.10) $\mathcal{G}_p(\eta) = E\#\mathcal{F}_p(\eta) = \{G : G = E\#F, F \in \mathcal{F}_p(\eta)\}, \qquad 0 < p < 2.$

In this paper, except where we explicitly state otherwise, the cdf's $F$ and $G$ are always related by scale mixing, so

$$G = E\#F.$$

(The relation $F \mapsto E\#F$ is one-to-one.) We often use $G$ and $G_n$ together, always implicitly assuming that they are related as the theoretical and empirical cdf of the same underlying samples so that $G_n$ is the empirical distribution for $n$ i.i.d. samples $X_i \sim G$, where

$$G_n(t) = \frac{1}{n} \sum_{i=1}^n \mathbb{1}_{\{X_i < t\}}.$$

**2. Asymptotics of minimax risk.** In this section, we prove Theorem 1.1. As usual, $R_n^*(M) = \sup_{\pi \in \Pi} \rho_n(\pi)$, where $\rho_n(\pi)$ denotes the Bayes risk $\mathcal{E}_\pi \mathcal{E}_\mu[\frac{1}{n}\|\log \hat{\mu}_\pi - \log \mu\|_2^2]$ with $\mu$ random, $\mu \sim \pi$; $\hat{\mu}_\pi$ denotes the Bayes estimator corresponding to prior $\pi$ and $\ell^2$ loss and $\Pi$ denotes the set of all priors supported on $M$ [here, $M = M_{n,p}(\eta)$, as in (1.3)]. Throughout this paper, we always implicitly assume that $P_{\pi_i}\{\mu_i \geq 1\} = 1$, where $\pi_i$ is the $i$th entry of $\pi$.

As in [7], we obtain a simple approximation of $R_n^*$ by considering a minimax-Bayes problem in which $\mu$ is a random vector that is only required to belong to $M$ *on average*. We define the minimax-Bayes risk as follows:

(2.1) $$\bar{R}_n^*(M_{p,n}(\eta)) = \inf_{\hat{\mu}} \sup_\pi \left\{ \mathcal{E}_\pi \mathcal{E}_\mu \left[ \frac{1}{n}\|\log \hat{\mu} - \log \mu\|_2^2 \right] : \mathcal{E}_\pi\left[\frac{1}{n}\sum_{i=1}^n \log^p \mu_i\right] \leq \eta^p \right\}.$$

Since a degenerate prior distribution concentrated at a single point $\mu \in M_{p,n}(\eta)$ trivially satisfies the moment constraint, the minimax-Bayes risk is an upper bound for the minimax risk, that is,

$$R_n^*(M_{n,p}(\eta)) \leq \bar{R}_n^*(M_{n,p}(\eta)). \tag{2.2}$$

In fact, for large $n$, we have asymptotic equality; in Section 2.1 we will prove the following:

THEOREM 2.1.
$$\lim_{n \to \infty} \frac{R_n^*(M_{n,p}(\eta))}{\bar{R}_n^*(M_{n,p}(\eta))} = 1.$$

Consider a univariate decision problem with data $X$ a scalar random variable, with $\mu$ a random scalar satisfying $\mu \sim F$ and $X|\mu \sim \text{Exp}(\mu)$. The corresponding *univariate* minimax-Bayes risk is

$$\bar{\rho}(\eta) = \bar{\rho}_p(\eta) = \inf_{\delta} \sup_{F \in \mathcal{F}_p(\eta)} \mathcal{E}_F \mathcal{E}_\mu (\log \delta(X) - \log \mu)^2. \tag{2.3}$$

The univariate and $n$-variate minimax risks are closely connected; in Section 2.2, we will prove the following:

THEOREM 2.2.   $\bar{R}_n^*(M_{n,p}(\eta)) = \bar{\rho}_p(\eta)$.

The univariate minimax-Bayes risk has a simple asymptotic expression as given by the following result:

THEOREM 2.3.  *For $0 < p < 2$,*
$$\lim_{\eta \to 0} \left( \frac{\bar{\rho}_p(\eta)}{\eta^p \log^{2-p} \log \frac{1}{\eta}} \right) = 1.$$

Theorem 1.1 follows immediately by combining Theorems 2.1–2.3.

2.1. *Proof of Theorem 2.1.*  Because (2.2) gives half of what we need, our task is to establish an asymptotic inequality in the other direction. We use a strategy similar to that of [7].

Now, for fixed $\eta$, choose $0 < \zeta \ll \eta$ and construct the product distribution $\Pi_{\eta-\zeta}^{(n)} = \prod_{i=1}^n \pi_{\eta-\zeta}^*$, where $\mu_i \overset{iid}{\sim} \pi_{\eta-\zeta}^*$, $\int \log^p(\mu) \, d\pi^* = (\eta - \zeta)^p$, $1 \leq i \leq n$, and $\pi^*$ is least favorable for univariate Bayes Minimax problem (2.3), so $\Pi_{\eta-\zeta}^{(n)}$ is least favorable for the $n$-variate Bayes Minimax problem (2.1). Let $A_n = \{\frac{1}{n} \sum_{i=1}^n \log^p \mu_i \leq \eta^p\}$. We then construct a new prior, $\tilde{\Pi}_{\eta-\zeta}^{(n)} = \Pi_{\eta-\zeta}^{(n)}(\cdot | A_n)$. By the law of large numbers (LLN),

$$P(A_n) \to 1, \tag{2.4}$$

while under $\Pi_{\eta-\zeta}^{(n)}$, we have $\mu \in M_{n,p}(\eta)$, that is, $\operatorname{supp} \Pi_{\eta-\zeta}^{(n)} \subset M_{n,p}(\eta)$. As the minimax risk is the supremum of Bayes risks, we have

(2.5) $$R_n^* \geq \rho_n(\tilde{\Pi}_{\eta-\zeta}^{(n)}).$$

Now, for any constant $w > 1$, and with $L(\cdot, \cdot)$ the loss function

$$L(\hat{\mu}, \mu) = \frac{1}{n} \sum_{i=1}^{n} (\log \hat{\mu}_i - \log \mu_i)^2,$$

define the $w$-truncated loss function,

$$L^{(w)}(\hat{\mu}, \mu) = \frac{1}{n} \sum_{i=1}^{n} \min\{(\log \hat{\mu}_i - \log \mu_i)^2, w\}.$$

Clearly,

(2.6) $$\rho_n(\tilde{\Pi}_{\eta-\zeta}^{(n)}, L) \geq \rho_n(\tilde{\Pi}_{\eta-\zeta}^{(n)}, L^{(w)}),$$

where $\rho_n(\pi, L)$ denotes the Bayes risk with respect to the loss function $L$. With $\|\cdot\|_{TV}$ denoting the variation distance, the definition of $\tilde{\Pi}_{\eta-\zeta}^{(n)}$ and (2.4) give

$$\|\tilde{\Pi}_{\eta-\zeta}^{(n)} - \Pi_{\eta-\zeta}^{(n)}\|_{TV} \leq 1 - P(A_n) \to 0.$$

For variation distance, $|\mathcal{E}_P f - \mathcal{E}_Q f| \leq \|f\|_\infty \cdot \|P - Q\|_{TV}$; thus, for any fixed $w$, the Bayes risk

$$|\rho_n(\tilde{\Pi}_{\eta-\zeta}^{(n)}, L^{(w)}) - \rho_n(\Pi_{\eta-\zeta}^{(n)}, L^{(w)})| \leq w \cdot (1 - P(A_n)) \to 0 \quad \text{as } n \to \infty.$$

On the other hand, for $L$ or $L^{(w)}$, the coordinatewise separability of the loss and the independence of the coordinates ensure that the per-coordinate Bayes risk does not depend on the number of coordinates, that is,

$$\rho_n(\Pi_{\eta-\zeta}^{(n)}, L) = \rho_1(\pi_{\eta-\zeta}^*, L), \qquad \rho_n(\Pi_{\eta-\zeta}^{(n)}, L^{(w)}) = \rho_1(\pi_{\eta-\zeta}^*, L^{(w)}).$$

We conclude that for each $w > 0$,

$$\rho_n(\tilde{\Pi}_{\eta-\zeta}^{(n)}, L^{(w)}) \to \rho_1(\pi_{\eta-\zeta}^*, L^{(w)}) \quad \text{as } n \to \infty.$$

Using monotone convergence of $L^{(w)} \to L$ as $w \to \infty$, we have

$$\rho_1(\pi_{\eta-\zeta}^*, L^{(w)}) \to \rho_1(\pi_{\eta-\zeta}^*, L) = \bar{\rho}(\eta - \zeta),$$

so from (2.5)–(2.6),

$$R_n^* \geq \bar{\rho}(\eta - \zeta).$$

Now, $\bar{\rho}(\eta)$ is monotone and continuous as a function of $\eta$; thus, by letting $\zeta \to 0$, we have

$$R_n^* \geq \bar{\rho}(\eta) = \bar{R}_n^*.$$

2.2. *Proof of Theorem 2.2.* First, observe that by the coordinatewise-separable nature of any estimator $\delta = \delta^n$ for $\mu$ and the i.i.d. structure of the $X_i/\mu_i$,

$$(2.7) \quad \frac{1}{n}\mathcal{E}_\pi \mathcal{E}_\mu \|\log \delta^n - \log \mu\|_2^2 = \frac{1}{n}\sum_i \int \mathcal{E}_{\mu_i}[\log \delta(X_i) - \log \mu_i]^2 \pi_i(d\mu_i)$$

$$(2.8) \qquad\qquad = \frac{1}{n}\int \mathcal{E}_{\mu_1}[\log \delta(X_1) - \log \mu_1]^2 \Big(\sum_i \pi_i\Big)(d\mu_1)$$

$$(2.9) \qquad\qquad = \mathcal{E}_{F_\pi}\mathcal{E}_{\mu_1}[\log \delta(X_1) - \log \mu_1]^2,$$

where $F_\pi = \frac{1}{n}\sum \pi_i(d\mu_1)$ is a univariate prior. Second, observe that the moment condition on $\pi$ can also be expressed in terms of $F_\pi$ since

$$(2.10) \quad \frac{1}{n}\mathcal{E}_\pi \sum \log^p \mu_i = \frac{1}{n}\sum_i \int \log^p(\mu_i)\pi_i(d\mu_i) = \int \log^p(\mu_1) F_\pi(d\mu_1),$$

thus, $\mathcal{E}_{F_\pi}\log^p \mu_1 \leq \eta^p$. Theorem 2.2 is easily derived from (2.7)–(2.10). Indeed, let $(F^0, \delta^0)$ be a saddlepoint for the univariate problem (2.3), that is, $\delta^0$ is a minimax rule, $F^0$ is a least favorable prior distribution and $\delta^0$ is Bayes for $F^0$. Let $F^{0,n}$ denote the $n$-fold Cartesian product measure derived from $F^0$ and $\delta^{0,n}$ the $n$-fold Cartesian product of $\delta^0$. From (2.10) and (2.7), $F^{0,n}$ satisfies the moment constraint for $\bar{R}_n^*(M_{n,p}(\eta))$ and

$$\frac{1}{n}\mathcal{E}_{F^{0,n}}\mathcal{E}_\mu \|\log \delta^{0,n} - \log \mu\|_2^2 = \bar{\rho}_p(\eta).$$

To establish the theorem, it is enough to verify that $(F^{0,n}, \delta^{0,n})$ is a saddlepoint for the minimax problem $\bar{R}_n^*(M_{n,p}(\eta))$. This would follow if for every $\pi$ obeying the moment constraint for $\bar{R}_n^*(M_{n,p}(\eta))$,

$$\mathcal{E}_\pi \mathcal{E}_\mu \|\log \delta^{0,n} - \log \mu\|_2^2 \leq \mathcal{E}_{F^{0,n}}\mathcal{E}_\mu \|\log \delta^{0,n} - \log \mu\|_2^2.$$

But (2.7)–(2.10) reduce this to the saddlepoint property of $(F^0, \delta^0)$ in the 1-dimensional minimax problem $\bar{\rho}_p(\eta)$.

2.3. *Proof of Theorem 2.3.* The following is proved in [11], Chapter 6:

LEMMA 2.1. *For functions $a = a(\eta)$ and $d = d(\eta)$ such that $\lim_{\eta \to 0} a(\eta) = 0$, $\lim_{\eta \to 0} d(\eta) = \infty$ and $\lim_{\eta \to 0}[a(\eta)/d(\eta)]^{1/(d(\eta)-1)} = 0$,*

$$\int_0^1 [(a/d) + y^{1-1/d}]^{-1}\,dy = d \cdot (1 + O((a/d)^{1/(d-1)})) \qquad as\ \eta \to 0.$$

We now describe lower and upper bounds for $\bar{\rho}(\eta)$, both equivalent to $\eta^p \log^{2-p}(\log\frac{1}{\eta})$ asymptotically as $\eta \to 0$. First, consider a lower bound for $\bar{\rho}(\eta)$. A natural lower bound uses 2-point priors,

$$(2.11) \qquad \bar{\rho}(\eta) \equiv \sup_{F \in \mathcal{F}_p(\eta)} \rho_1(F) \geq \sup_{\{(\varepsilon,\mu):\, \varepsilon \log^p(\mu) = \eta^p\}} \rho_1(F_{\varepsilon,\mu}),$$

where $F_{\varepsilon,\mu} = (1-\varepsilon)\nu_1 + \varepsilon \nu_\mu \in \mathcal{F}_p(\eta)$ denotes the mixture of mixing point masses at 1 and $\mu$ with fractions $(1-\varepsilon)$ and $\varepsilon$, respectively. The Bayes rule $\delta_B(X; F_{\varepsilon,\mu})$ obeys

$$(2.12) \qquad \log(\delta_B(X; F_{\varepsilon,\mu})) = \frac{\frac{\varepsilon}{\mu} e^{-X/\mu}}{(1-\varepsilon)e^{-X} + \frac{\varepsilon}{\mu}e^{-X/\mu}} \log \mu$$

and the Bayes risk is

$$\rho_1(F_{\varepsilon,\mu}) = (\log \mu)^2 \int_0^\infty \frac{(1-\varepsilon)e^{-x}\frac{\varepsilon}{\mu}e^{-\frac{x}{\mu}}}{(1-\varepsilon)e^{-x} + \frac{\varepsilon}{\mu}e^{-\frac{x}{\mu}}}\, dx$$

$$(2.13) \qquad = \frac{\varepsilon \log^2(\mu)}{\mu} \int_0^1 \left(\frac{\varepsilon}{(1-\varepsilon)\mu} + y^{1-\frac{1}{\mu}}\right)^{-1} dy;$$

particularly, if we let $\mu^* = \mu^*(\eta) = \log(\frac{1}{\eta})/(\log\log\frac{1}{\eta})$ and $\varepsilon^* = \varepsilon^*(\eta) = \eta^p/\log^p(\mu^*)$, then applying Lemma 2.1 with $a = a(\eta) = \varepsilon^*/(1-\varepsilon^*)$ and $d = d(\eta) = \mu^*$, we have

$$\rho_1(F_{\varepsilon^*(\eta),\mu^*(\eta)}) = \left(\eta^p \log^{2-p} \log \frac{1}{\eta}\right) \cdot (1 + o(1))$$

and obtain the desired lower bound

$$(2.14) \qquad \bar{\rho}(\eta) \geq \rho_1(F_{\varepsilon^*(\eta),\mu^*(\eta)}) = \left(\eta^p \log^{2-p} \log \frac{1}{\eta}\right) \cdot (1 + o(1)).$$

We obtain an upper bound by considering the risk of thresholding. Define the univariate thresholding nonlinearity

$$(2.15) \qquad \delta_t(x) = \begin{cases} x, & x \geq t, \\ 1, & \text{otherwise.} \end{cases}$$

Then with thresholding estimator $\delta_t(X)$ based on scalar data $X$ obeying $X|\mu \sim \text{Exp}(\mu)$, where the scalar $\mu$ is distributed according to a prior $F \in \mathcal{F}_p(\eta)$, the univariate Bayes thresholding risk is

$$\rho_T(t, F) = \mathcal{E}(\log(\delta_t(X)) - \log(\mu))^2.$$

We are particularly interested in the specific threshold

$$t_0 = t_0(p, \eta) = p\log\left(\frac{1}{\eta}\right) + p\log\log\left(\frac{1}{\eta}\right) + \sqrt{\log\log\left(\frac{1}{\eta}\right)}.$$

The worst-case univariate Bayes risk for this rule is

(2.16) $$\bar{\rho}_T(t_0, \eta) = \bar{\rho}(t_0, \eta; p) \equiv \sup_{F \in \mathcal{F}_p(\eta)} \rho_T(t_0, F).$$

As the minimax rule is at least as good as any specific rule, we have

(2.17) $$\bar{\rho}(\eta) \leq \bar{\rho}_T(t_0, \eta).$$

Now, in the proof of Theorem 1.2 below, we show that the thresholding risk obeys

(2.18) $$\bar{\rho}_T(t_0, \eta; p) \leq \eta^p \log^{2-p} \frac{1}{\eta}(1 + o(1)) \qquad \text{as } \eta \to 0.$$

Combining the lower bound given by (2.14) and the upper bounds given by (2.17)–(2.18), we obtain Theorem 2.3.

**3. Asymptotic minimaxity of thresholding.** We now prove Theorem 1.2, showing that thresholding estimates can asymptotically approach the minimax risk.

3.1. *Reduction to univariate thresholding.* In effect, we need only prove (2.18). We first recall why this establishes Theorem 1.2. Again, let $\hat{\mu}_t$ denote the thresholding procedure on samples of size $n$. Trivially, for any $t$ and $n$, the risk of thresholding at $t$ exceeds the minimax risk

$$\sup_{M_{n,p}(\eta)} R_n(\hat{\mu}_t, \mu) \geq R_n^*(M_{n,p}(\eta)).$$

Theorem 1.2 thus follows from an asymptotic inequality in the other direction,

(3.1) $$\limsup_{\eta \to 0} \inf_t \left[ \limsup_{n \to \infty} \frac{\sup_{M_{n,p}(\eta)} R_n(\hat{\mu}_t, \mu)}{R_n^*(M_{n,p}(\eta))} \right] \leq 1.$$

If we take

(3.2) $$t_0 = t_0(p, \eta) = p \log(1/\eta) + p \log \log(1/\eta) + \sqrt{\log \log(1/\eta)},$$

then by Theorem 2.1 and Theorem 2.2, (3.1) reduces to

(3.3) $$\limsup_{\eta \to 0} \left[ \frac{\limsup_{n \to \infty} \sup_{M_{n,p}(\eta)} R_n(\hat{\mu}_{t_0}, \mu)}{\bar{\rho}(\eta)} \right] \leq 1.$$

Consider the worst Bayes risk of $\hat{\mu}_{t_0}$ with respect to any prior $\mu \sim \pi$, where $\pi$ is the distribution of a random vector which is only required to belong to $M_{n,p}$ on average,

$$\bar{R}_n^*(\hat{\mu}_{t_0}, \eta) = \bar{R}_n^*(\hat{\mu}_{t_0}, \eta; p)$$

(3.4)
$$= \sup \left\{ \mathcal{E}_\pi \mathcal{E}_\mu \left[ \frac{1}{n} \|\log \hat{\mu}_{t_0} - \log \mu\|_2^2 \right], \text{ for } \pi : \mathcal{E}_\pi \frac{1}{n} \sum_{i=1}^n \log^p \mu_i \leq \eta^p \right\}.$$

Now, since degenerate prior distributions concentrated at points $\mu \in M_{p,n}(\eta)$ trivially satisfy the moment constraint $\mathcal{F}_p(\eta)$, we have

$$\sup_{M_{n,p}(\eta)} R_n(\hat{\mu}_{t_0}, \mu) \leq \bar{R}_n^*(\hat{\mu}_{t_0}, \eta). \tag{3.5}$$

Consider also the worst univariate Bayes risk (2.16) of the scalar rule $\delta_{t_0}(X)$, as in (2.15), with respect to univariate prior $F \in \mathcal{F}_p(\eta)$. As in the proof of Theorem 2.2, it is not hard to show that the minimax multivariate Bayes risk is the same as the minimax univariate Bayes risk

$$\bar{R}_n^*(\hat{\mu}_{t_0}, \eta) = \bar{\rho}_T(t_0, \eta). \tag{3.6}$$

Hence, we now see that given (2.14), the matching upper bound (2.18) implies that

$$\lim_{\eta \to 0} \frac{\bar{\rho}_T(t_0, \eta)}{\bar{\rho}(\eta)} = 1. \tag{3.7}$$

Combining (3.5)–(3.7) yields (3.3) and Theorem 1.2. We thus turn to (2.18).

The univariate Bayes risk for thresholding at $t$ can be decomposed into a *bias proxy* and a *variance proxy* as follows:

$$\bar{\rho}_T(t, F) = \int (\log \mu)^2 (1 - e^{-\frac{t}{\mu}}) \, dF(\mu) + \int \left[ \int_{\frac{t}{\mu}}^{\infty} \log^2(x) e^{-x} \, dx \right] dF(\mu),$$

$$\equiv \int b(t, \mu) \, dF(\mu) + \int v(t, \mu) \, dF(\mu),$$

say. We now proceed to show that as $\eta \to 0$,

$$\sup_{F \in \mathcal{F}_p(\eta)} \int b(t_0, \mu) \, dF(\mu) \leq \eta^p \log^{2-p} \log \frac{1}{\eta} \tag{3.8}$$

and

$$\sup_{F \in \mathcal{F}_p(\eta)} \int v(t_0, \mu) \, dF(\mu) = o\left(\eta^p \log^{2-p} \log \frac{1}{\eta}\right). \tag{3.9}$$

Together, these imply (2.18).

3.2. *Maximizing linear functionals over* $\mathcal{F}_p(\eta)$. The relations (3.8)–(3.9) concern maximization of functionals over cdf's of moment-constrained scale mixtures. We now approach this problem from a general viewpoint, looking ahead to maximization problems in later sections.

Consider two functions $\psi(\mu), \phi(\mu)$ in $C[1, \infty) \cap C^2(1, \infty)$. Suppose

(a) $\phi$ is strictly increasing and $\phi(1) = 0$;
(b) $\psi$ is bounded, $\psi(1) = 0$, $\psi \geq 0$ but $\psi$ is not identically 0;
(c) $\lim_{\mu \to \infty} [\psi(\mu)/\phi(\mu)] = 0$.

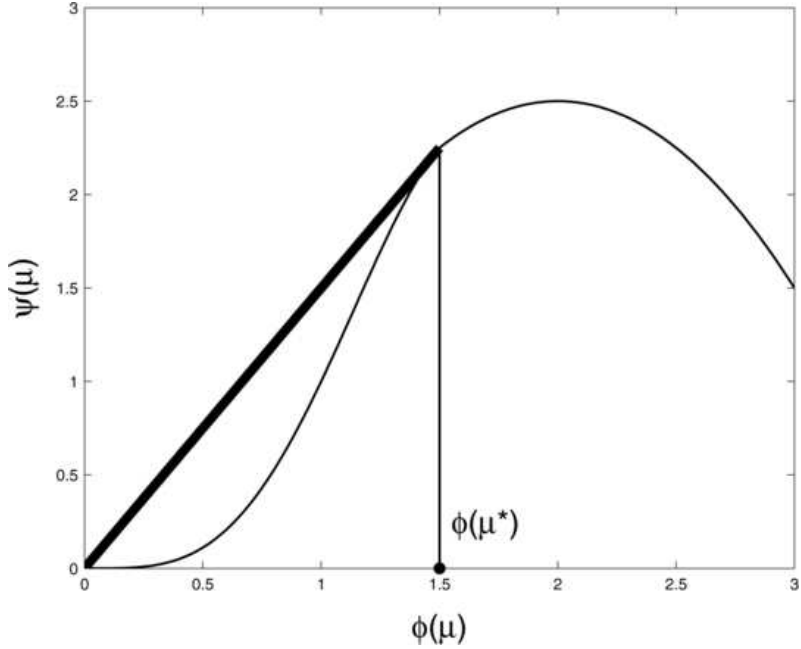

FIG. 1. *Generalized convex envelope $\Psi(z)$ for the case $\lim_{\mu \to 1+}[\psi(\mu)/\phi(\mu)] < \infty$ in the $\phi$–$\psi$ plane. In this example shown here with $\lim_{\mu \to 1+}[\psi(\mu)/\phi(\mu)] = 0$, the thinner curve is $\{(\phi(\mu), \psi(\mu)) : \mu \geq 1\}$. When $0 \leq z \leq \phi(\mu^*)$, $\Psi(z)$ is a linear function of $z$ and is illustrated by the line segment. The case $z > \phi(\mu^*)$ is not discussed.*

We are interested in the following maximization problem:

$$\text{(3.10)} \quad \Psi(z) = \sup_{F \in \mathcal{F}} \left\{ \int \psi(\mu) \, dF(\mu) : \int \phi(\mu) \, dF(\mu) \leq z \right\}.$$

In the case $\phi(\mu) = \mu$, $\Psi(z)$ is the usual convex envelope of $\psi$, that is, $\Psi(z)$ traces out the least concave majorant of the graph of $\Psi$. The next two lemmas describe the computation of the envelope.

LEMMA 3.1. *Suppose $\lim_{\mu \to 1+}[\psi(\mu)/\phi(\mu)]$ exists and the limit is strictly smaller than $\Psi^* \equiv \sup_{\mu > 1}\{\psi(\mu)/\phi(\mu)\}$. Set*

$$\mu^* = \mu^*(\psi, \phi) \equiv \max\{\mu > 1 : \psi(\mu)/\phi(\mu) = \Psi^*\}.$$

*Then for any $0 \leq z \leq \phi(\mu^*)$, $\Psi(z) = \Psi^* \cdot z$ and is attained by the mixture of point masses at $1$ and $\mu^*$ with masses $(1 - \varepsilon(z))$ and $\varepsilon(z)$, respectively, where $\varepsilon(z) = \varepsilon(z; \psi, \phi) = z/\phi(\mu^*)$.*

See Figure 1.

LEMMA 3.2. *Suppose that* $\lim_{\mu \to 1+}[\psi(\mu)/\phi(\mu)] = \infty$ *and suppose there exists* $\bar{\mu} = \bar{\mu}(\psi, \phi) > 1$ *so that* $(\psi'(\mu)/\phi'(\mu))$ *is strictly decreasing in the interval* $(1, \bar{\mu}]$ *and, finally, that* $\psi'(\bar{\mu})/\phi'(\bar{\mu}) < \Psi^{**}(\bar{\mu})$, *where*

$$(3.11) \quad \Psi^{**}(\mu) = \Psi^{**}(\mu; \bar{\mu}, \phi, \psi) \equiv \sup_{\mu' > \bar{\mu}} \frac{\psi(\mu') - \psi(\mu)}{\phi(\mu') - \phi(\mu)}, \qquad 1 \leq \mu < \bar{\mu}.$$

*Then there is a unique solution* $\mu_* = \mu_*(\psi, \phi)$ *to the equation*

$$\Psi^{**}(\mu) = \psi'(\mu)/\phi'(\mu), \qquad 1 < \mu \leq \bar{\mu};$$

*moreover, letting*

$$\mu^* = \max\left\{\mu \geq \bar{\mu} : \frac{\psi(\mu) - \psi(\mu_*)}{\phi(\mu) - \phi(\mu_*)} = \Psi^{**}(\mu_*)\right\},$$

*then when* $0 < z \leq \phi(\mu_*)$, $\Psi(z) = \psi(\phi^{-1}(z))$ *and is attained by the single point mass* $\nu_{\mu_z}$ *with* $\mu_z = \phi^{-1}(z)$ *and when* $\phi(\mu_*) < z \leq \phi(\mu^*)$, $\Psi(z) = \psi(\mu_*) + \Psi^{**}(\mu_*)[z - \phi(\mu_*)]$ *and is attained by the mixture of point masses at* $\mu_*$ *and* $\mu^*$ *with masses* $(1 - \varepsilon(z))$ *and* $\varepsilon(z)$, *respectively, where* $\varepsilon(z) = \varepsilon(z; \phi, \psi) = [z - \phi(\mu_*)]/[\phi(\mu^*) - \phi(\mu_*)]$.

Notice here that the strict monotonicity of $\psi'(\mu)/\phi'(\mu)$ over $(1, \bar{\mu}]$ is equivalent to concavity of the curve $\{(\phi(\mu), \psi(\mu)) : 1 < \mu \leq \bar{\mu}\}$ in the $(\phi(\mu), \psi(\mu))$ plane. See Figure 2.

The proofs of Lemmas 3.1 and 3.2 can be found in the full version of this paper [6].

3.3. *Maximizing bias and variance.* To apply Lemma 3.1 to the bias proxy, set $\psi = \psi_\eta(\mu) = b(t_0, \mu) = \log^2(\mu)(1 - e^{-t_0/\mu})$, $\phi(\mu) = \log^p(\mu)$ and $\Psi(z)$, as in (3.10). Then the worst bias $\sup_{\mathcal{F}_p(\eta)} \int b(t_0, \mu) \, dF \equiv \Psi(\eta^p)$. Direct calculation shows that for large $t_0$,

$$\mu^* \equiv \operatorname{argmax}[\psi(\mu)/\phi(\mu)] \sim \frac{t_0}{\log \log t_0 - \log(2 - p)}$$

and

$$\Psi^* = \bar{\Psi}_{p,\eta} \equiv \frac{\psi(\mu^*)}{\log^p(\mu^*)} \sim \log^{2-p} t_0 \sim \log^{2-p}\left(\frac{1}{\eta}\right).$$

It is obvious that for sufficiently small $\eta$, $\eta^p < \phi(\mu^*)$; thus, by Lemma 3.1, $\Psi(\eta^p) = \Psi^* \cdot \eta^p$ and relation (3.8) follows directly.

Now consider the variance proxy. Let $\psi(\mu) = \psi_\eta(\mu) \equiv v(t_0, \mu) - v(t_0, 1)$, $\phi(\mu) = \log^p(\mu)$ and again with $\Psi(z)$ as in (3.10), the maximal variance proxy $\sup_{\mathcal{F}_p(\eta)} \int v(t_0, \mu) \, dF = \Psi(\eta^p) + v(t_0, 1)$. Notice here that $v(t_0, 1) = o(\eta^p \log^{2-p}(\log \frac{1}{\eta}))$, so to show relation (3.9), we need only demonstrate that

$$(3.12) \qquad \qquad \Psi(\eta^p) = O(\eta^p).$$

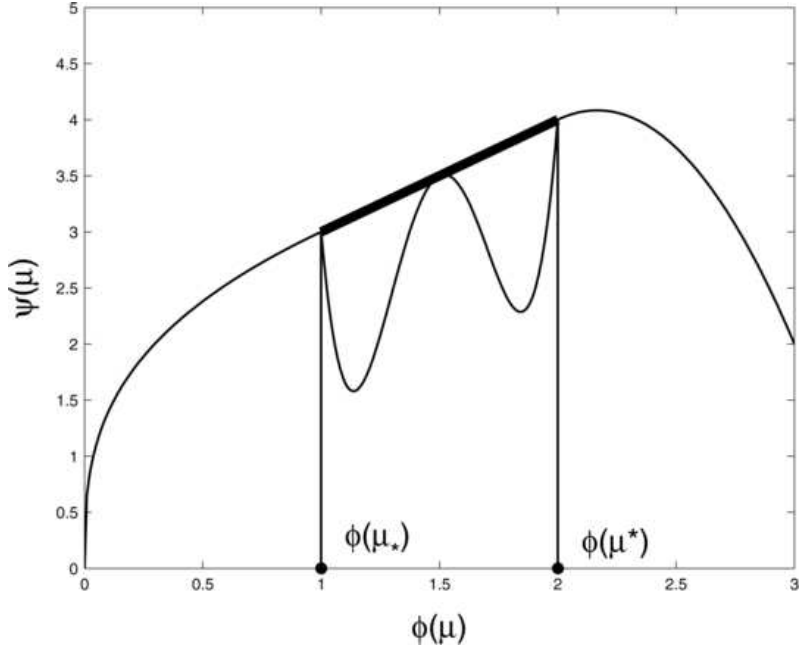

Fig. 2. *Generalized convex envelope $\Psi(z)$ for the case $\lim_{\mu \to 1+}[\psi(\mu)/\phi(\mu)] = \infty$ in the $\phi$-$\psi$ plane. The thinner curve is $\{(\phi(\mu), \psi(\mu)) : \mu \geq 1\}$. When $0 < \mu < \mu_*$, $\{(\phi(\mu), \Psi(\mu)) : 0 < \mu < \mu_*\}$ traces out the same curve as that of $\{(\phi(\mu), \psi(\mu)) : 0 < \mu < \mu_*\}$ and when $\mu_* \leq \mu \leq \mu^*$, $\Psi(z)$ is a linear function of $z = \phi(\mu)$ which is illustrated by the line segment. The slope of the line segment equals that of the tangent at $\mu_*$ of the curve $\{(\phi(\mu), \psi(\mu)) : \mu \geq 1\}$. The case $z > \phi(\mu^*)$ is not discussed.*

Direct calculations show that

$$(3.13) \qquad \lim_{\mu \to 1+}\left[\frac{\psi(\mu)}{\phi(\mu)}\right] = \begin{cases} 0, & 0 < p < 1, \\ t_0 \log^2(t_0)e^{-t_0}, & p = 1, \\ \infty, & 1 < p < 2, \end{cases}$$

so we will calculate $\Psi(z)$ for the cases $0 < p \leq 1$ and $1 < p < 2$ separately.

When $0 < p \leq 1$, let $c = \int_1^\infty \log^2(x) e^{-x}\, dx$ and note that for sufficiently large $t_0$, the condition of Lemma 3.1 is satisfied; moreover, direct calculations show that

$$\mu^* = \underset{\mu > 1}{\operatorname{argmax}}\{\psi(\mu)/\phi(\mu)\} \sim t_0, \qquad \Psi^* = \frac{\psi(\mu^*)}{\log^p(\mu^*)} \sim \frac{c}{\log^p(t_0)};$$

for sufficiently small $\eta$, we have $\eta^p < \phi(\mu^*)$, so by Lemma 3.1, $\Psi(\eta^p) = \Psi^* \cdot \eta^p$ and (3.12) follows directly.

When $1 < p < 2$, if we let $\bar{\mu}$ denote the smaller solution of the equation $\frac{t_0}{\mu}\log(\mu) = (p-1)$, then for large $t_0$, $\bar{\mu} \sim 1 + \frac{p-1}{t_0}$; moreover, by elementary analysis, $[\psi'(\mu)/\phi'(\mu)]$ is strictly decreasing in $(1, \bar{\mu}]$ and $\psi'(\bar{\mu})/\phi'(\bar{\mu}) <$

$\Psi^{**}(\bar{\mu})$ and the condition of Lemma 3.2 is satisfied. Furthermore, for large $t_0$,

$$\Psi^{**}(\mu) \sim \frac{c}{\log^p t_0}, \qquad \forall 1 < \mu \leq \bar{\mu}. \tag{3.14}$$

More elementary analysis shows that

$$\mu^* = \operatorname*{argmax}_{\mu \geq \bar{\mu}} \frac{\psi(\mu) - \psi(\mu_*)}{\phi(\mu) - \phi(\mu_*)} \sim \operatorname*{argmax}_{\mu \geq \bar{\mu}} \frac{\psi(\mu)}{\phi(\mu)} \sim t_0$$

and

$$\mu_* = \exp([ct_0 \log^{2+p} t_0 e^{-t_0}/p]^{1/(p-1)}),$$

$$\phi(\mu_*) = [ct_0 \log^{2+p} t_0 e^{-t_0}/p]^{p/(p-1)}.$$

It is now clear that for sufficiently small $\eta > 0$, $\phi(\mu_*) < \eta^p < \phi(\mu^*)$. Thus, by Lemma 3.2,

$$\Psi(\eta^p) = \psi(\mu_*) + \Psi^{**}(\mu_*)(\eta^p - \log(\mu_*)). \tag{3.15}$$

Taking $\mu = \mu_*$ in (3.14) and (3.15) gives (3.12), since

$$\Psi(\eta^p) = \psi(\mu_*) + \Psi^{**}(\mu_*)[\eta^p - \phi(\mu_*)] \sim \eta^p \frac{c}{\log^p t_0} = o(\eta^p).$$

**4. The FDR functional.** We now come to the central idea in our analysis of FDR thresholding—to view the FDR threshold as a functional of the underlying cumulative distribution function. For any fixed $0 < q < 1$, the *FDR functional* $T_q(\cdot)$ is defined as

$$T_q(G) = \inf\left\{t : \bar{G}(t) \geq \frac{1}{q}\bar{E}(t)\right\}, \tag{4.1}$$

where $G$ is any cdf.

The relevance of $T_q$ follows from a simple observation. If $G_n$ is the empirical distribution of $X_1, X_2, \ldots, X_n$, then $T_q(G_n)$ is effectively the same as the FDR threshold $\hat{t}^{FDR}(X_1, \ldots, X_n)$. More precisely (see Lemma 6.1 below), thresholding at $T_q(G_n)$ and at $\hat{t}^{FDR}(X_1, \ldots, X_n)$ always gives, numerically, exactly the same estimate $\hat{\mu}_{q,n}$.

In this section, we consider several key properties of this functional.

4.1. *Definition, boundedness and continuity.* We first observe that $T_q(G)$ is well defined at nontrivial scale mixtures of exponentials.

LEMMA 4.1 (Uniqueness). *For fixed $0 < q < 1$ and for all $G \in \mathcal{G}$, $G \neq E$, the equation*

$$\bar{G}(t) = \frac{1}{q}\bar{E}(t) \tag{4.2}$$

*has a unique solution on $[0, \infty)$ which we denote $T_q(G)$.*

PROOF. Indeed, with $\mu$ a random variable greater than or equal to 1, $\bar{G}(t) = \mathcal{E}[\bar{E}(t/\mu)]$. Hence, if $\mu \neq 1$ a.s., then for some $\mu_0 > 1$ and some $\varepsilon > 0$, we have that for all $t \geq 0$, $\bar{G}(t) > \varepsilon \bar{E}(t/\mu_0)$. Now, $\bar{G}(0) < \bar{E}(0)/q$, while for sufficiently large $t$, $\bar{E}(t)/q < \varepsilon \bar{E}(t/\mu_0)$. Hence, for some $t = t_0$ on $[0, \infty)$, (4.2) holds. Now, consider the slope of $\bar{G}(t)$,

$$-\frac{d}{dt}\bar{G}(t) = \mathcal{E}[\bar{E}(t/\mu)/\mu] < \mathcal{E}[\bar{E}(t/\mu)] = \bar{G}(t).$$

Compare this with the slope of $\bar{E}(t)/q$. We have

$$-\frac{d}{dt}\frac{1}{q}\bar{E}(t) = \frac{1}{q}\bar{E}(t).$$

At $t = t_0$, $\frac{1}{q}\bar{E}(t) = \bar{G}(t)$, so

$$\frac{d}{dt}\left(\bar{G}(t_0) - \frac{1}{q}\bar{E}(t)\right)\bigg|_{t=t_0} > 0.$$

In short, at any crossing of $\bar{G} - \frac{1}{q}\bar{E}$, the slope is positive. Downcrossings being impossible, there is only one upcrossing, so the solution (4.2) is unique. □

The ideas used in the proof immediately lead to two other important properties of $T_q$.

LEMMA 4.2 (Quasi-Concavity). *The collection of distributions $G \in \mathcal{G}$ satisfying $T_q(G) = t$ is convex. The collection of distributions satisfying $T_q(G) \geq t$ is convex.*

PROOF. The uniqueness lemma shows that the set $T_q(G) = t$ consists precisely of those cdf's $G$ obeying $\bar{G}(t) = e^{-t}/q$; this is a linear equality constraint over the convex set $\mathcal{G}$ and defines a convex subset of $\mathcal{G}$. The set $T_q(G) \geq t$ consists precisely of those cdf's $G$ obeying $\bar{G}(t) \leq e^{-t}/q$; this is a linear inequality constraint over the convex set $\mathcal{G}$ and generates a convex subset. □

We also immediately have the following:

LEMMA 4.3 (Stochastic Ordering). *We introduce the following notation for cdf's: $G_0 \lesssim G_1$ if $\bar{G}_1(t) \geq \bar{G}_0(t)$ for all $t > 0$. Then*

$$G_0 \lesssim G_1 \implies T_q(G_0) \geq T_q(G_1).$$

We now turn to boundedness and continuity of $T_q$. Recall that the Kolmogorov–Smirnov distance between cdf's $G$ and $G'$ is defined by

$$\|G - G'\| = \sup_t |G(t) - G'(t)|.$$

Viewing the collection of cdf's as a convex set in a Banach space equipped with this metric, the FDR functional $T_q(\cdot)$ is, in fact, locally bounded over neighborhoods of nontrivial scale mixture of exponentials.

LEMMA 4.4 (Boundedness). *For $G \in \mathcal{G}$, $G \neq E$,*

$$-\log\left(\frac{q}{1-q}\|G-E\|\right) \leq T_q(G) \leq \frac{1-q}{q}\frac{1}{\|G-E\|}.$$

PROOF. We introduce the shorthand notation $\tau = T_q(G)$. The left-hand inequality follows from $\bar{G}(\tau) = \bar{E}(\tau)/q$, which gives

$$\|G - E\| = \sup_t |G(t) - E(t)| \geq \bar{G}(\tau) - e^{-\tau} = \frac{1-q}{q}e^{-\tau}.$$

For the right-hand inequality, again use $\bar{G}(\tau) = \bar{E}(\tau)/q$ and convexity of $e^t$ to obtain

$$\frac{1}{q} = \int e^{(1-\frac{1}{\mu})\tau} dF \geq 1 + \tau \cdot \int \left(1 - \frac{1}{\mu}\right) dF.$$

At the same time, since $E \lesssim G$, we have $\|G - E\| = \sup_{t>0} \int [e^{-\frac{t}{\mu}} - e^{-t}] dF$. Observe that as a function of $t$, $\int [e^{-\frac{t}{\mu}} - e^{-t}] dF$ has a unique maximum point $t = \bar{t}$ satisfying $\int \frac{1}{\mu} e^{-\frac{\bar{t}}{\mu}} dF = e^{-\bar{t}}$, so

$$\|G - E\| = \int [e^{-\frac{\bar{t}}{\mu}} - e^{-\bar{t}}] dF = \int \left(1 - \frac{1}{\mu}\right) e^{-\frac{\bar{t}}{\mu}} dF \leq \int \left(1 - \frac{1}{\mu}\right) dF$$

and we have $\tau \leq \frac{1-q}{q} \frac{1}{\|G-E\|}$. □

In fact, the FDR functional is even locally Lipschitz away from $G = E$. Note that the image of the mapping $T_q : \mathcal{G} \mapsto \mathbb{R}$ is the interval $(\log(\frac{1}{q}), \infty)$.

LEMMA 4.5 (Modulus of Continuity). *Define*

$$\omega^*(\varepsilon; t_0) \equiv \sup\{|T_q(G') - t_0| : T_q(G) = t_0, \|G - G'\| \leq \varepsilon, G \in \mathcal{G}\}.$$

*Then for each fixed $t_0 > \log(1/q)$,*

(4.3) $$\omega^*(\varepsilon; t_0) \leq \frac{q}{\log(1/q)} t_0 e^{t_0} \varepsilon \cdot (1 + o(1)) \qquad \text{as } \varepsilon \to 0.$$

Crucially, the estimate (4.3) is uniform over $\{G \in \mathcal{G}, T_q(G) \leq t_0\}$ for fixed $t_0 > 0$. The proof even shows that

(4.4) $$\omega^*(\varepsilon; t_0) \leq C \cdot \varepsilon \qquad \text{for } 0 < \varepsilon < \varepsilon_{t_0},$$

where $C = C_{t_0,q} < \infty$ if $t_0 < \infty$. This implies the local Lipschitz property.

PROOF. Consider the optimization problem of finding the cdf $G^* \in \mathcal{G}$ which satisfies $T_q(G^*) = t_0$ and, subject to that constraint, is as 'steep' as possible at $t_0$, that is,

$$(4.5) \qquad \frac{\partial}{\partial t}\bar{G}^*(t)\bigg|_{t=t_0} = \inf\left\{\frac{\partial}{\partial t}\bar{G}(t)\bigg|_{t=t_0} : \bar{G}(t_0) = \frac{1}{q}\bar{E}(t_0), G \in \mathcal{G}\right\}.$$

Letting $\phi(\mu) = e^{-t_0/\mu}$ and $\psi(\mu) = (t_0/\mu)e^{-t_0/\mu}$, Problem (4.5) can be viewed as maximizing the linear functional $\int \psi(\mu)\,dF(\mu)$ with the constraint $\int \phi(\mu)\,dF(\mu) = \frac{1}{q}e^{-t_0}$. Observe that $\psi'(\mu)/\phi'(\mu)$ strictly decreases in $\mu$ over $(1, \infty)$, so in the $\phi$–$\psi$ plane, the curve $(\phi(\mu), \psi(\mu))$ is strictly concave and by arguments used in the proof of Lemma 3.2, the constrained maximum of $\int \psi(\mu)\,dF(\mu)$ is obtained at the point mass $F$ which satisfies $\int \phi(\mu)\,dF(\mu) = \frac{1}{q}e^{-t_0}$.

It thus follows that the solution to Problem (4.5) is $\bar{G}^*_{t_0}(t) = e^{-t/\mu^*}$ for $\mu^* = 1/(1 + \log(q)/t_0)$. It has the remarkable property that if $T_q(G) = t_0$,

$$(4.6) \qquad \bar{G}(t) \leq \bar{G}^*_{t_0}(t), \quad 0 < t < t_0, \qquad \bar{G}(t) \geq \bar{G}^*_{t_0}(t), \quad t > t_0.$$

Indeed, letting

$$h(t) \equiv [\bar{G}(t)/\bar{G}^*_{t_0}(t)] - 1 = \int e^{(\frac{1}{\mu^*} - \frac{1}{\mu})t}\,dF(\mu) - 1,$$

direct calculation shows that $h(t)$ is strictly convex as long as $P_F\{\mu = \mu^*\} \neq 1$ (otherwise $h \equiv 0$) and (4.6) follows by observing that $h(0) = h(t_0) = 0$.

For sufficiently small $\varepsilon$, define $t_-$ by

$$(4.7) \qquad \bar{G}^*_{t_0}(t_-) + \varepsilon = \bar{E}(t_-)/q$$

and define $t_+$ to be the smallest solution to the equation

$$(4.8) \qquad \bar{G}^*_{t_0}(t) - \varepsilon = \bar{E}(t)/q;$$

see Figure 3. Now, if $\|G' - G\| \leq \varepsilon$, then by (4.6) and (4.8),

$$\bar{G}'(t_+) \geq \bar{G}(t_+) - \varepsilon \geq \bar{G}^*_{t_0}(t_+) - \varepsilon = \bar{E}(t_+)/q,$$

hence, $T_q(G') \leq t_+$. Similarly, by (4.6) and (4.7),

$$(4.9) \qquad \bar{G}'(t_-) \leq \bar{G}(t_-) + \varepsilon \leq \bar{G}^*_{t_0}(t_-) + \varepsilon = \bar{E}(t_-)/q.$$

Observe that the function $(\bar{G}^*_{t_0}(t) - \bar{E}(t)/q)$ is strictly decreasing in the interval $[0, t_0]$, so (4.9) can be strengthened into

$$\bar{G}'(t) \leq \bar{G}(t) + \varepsilon \leq \bar{G}^*_{t_0}(t) + \varepsilon < \bar{E}(t)/q, \qquad 0 < t < t_-,$$

hence, $T_q(G') \geq t_-$. It follows that

$$(4.10) \qquad \omega(\varepsilon; t_0) \leq \max\{t_0 - t_-(\varepsilon), t_+(\varepsilon) - t_0\}.$$

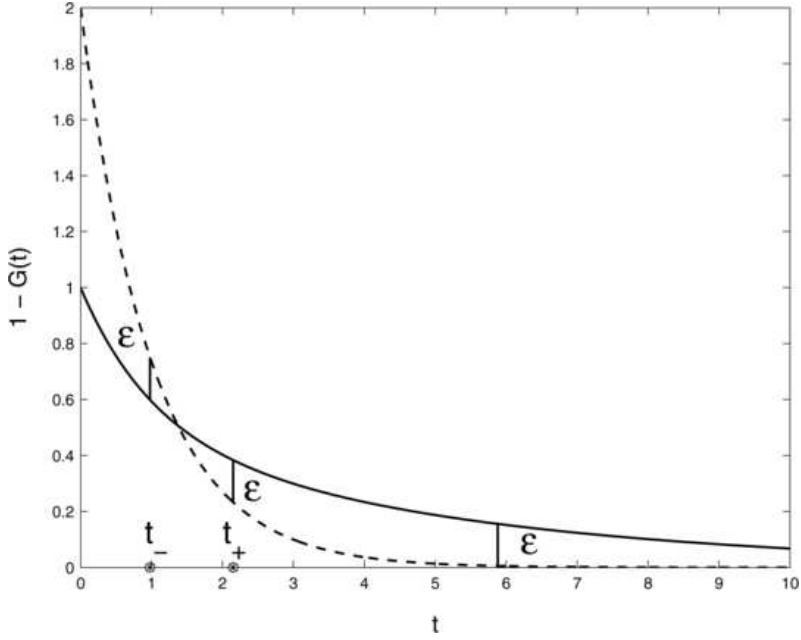

FIG. 3.  *The dashed curve is $(1/q)\bar{E}(t)$ with $q = 1/2$ and the solid curve is $\bar{G}^*_{t_0}(t)$. In the plot, $t_-$ is the solution of $\bar{G}^*_{t_0}(t) + \varepsilon = (1/q)\bar{E}(t)$ and $t_+$ is the smallest solution to the equation $\bar{G}^*_{t_0}(t) - \varepsilon = (1/q)\bar{E}(t)$. For any other $G$ with $T_q(G) = t_0$, $\bar{G}(t)$ is bounded above by $\bar{G}^*_{t_0}(t)$ when $0 < t < t_0$ and is bounded below by $\bar{G}^*_{t_0}(t)$ when $t > t_0$; moreover, for any $G'$ with $\|G' - G\| \leq \varepsilon$, $t_- \leq T_q(G') \leq t_+$.*

Finally, setting $w = t_+ - t_0$, (4.7) can be rewritten as $e^{-w/\mu^*} - e^{-w} = \varepsilon q e^{t_0}$. Letting $w(\delta)$ denote the smaller of the two solutions to $e^{-w/\mu^*} - e^{-w} = \delta$, elementary analysis shows that for small $\delta > 0$, $w(\delta) \sim \delta/(1 - 1/\mu^*) = \delta t_0 / \log(1/q)$, so as $\varepsilon \to 0$, $t_+ - t_0 \sim (q/\log(1/q)) \cdot t_0 e^{t_0} \varepsilon$ and, similarly, $t_0 - t_-(\varepsilon) \sim (q/\log(1/q)) \cdot t_0 e^{t_0} \varepsilon$. Inserting these into (4.10) gives the lemma.  □

4.2. *Behavior under the Bayesian model.* The continuity of $T_q$ established in Lemma 4.5 and the role of minimax Bayes risk in solving for the minimax risk in Sections 2 and 3 combine to suggest a fruitful change of viewpoint. Instead of viewing the $X_i \sim \text{Exp}(\mu_i)$ with fixed constants $\mu_i$, $i = 1, \ldots, n$, we view the $\mu_i$ as themselves sampled i.i.d. from a distribution $F$, so the $X_i$ are sampled i.i.d. from a mixture of exponentials $G = E \# F$. Starting now and continuing through Sections 5 and 6, we adopt this viewpoint exclusively. Moreover, for our sparsity constraint, instead of assuming that $\frac{1}{n}(\sum_{i=1}^{n}(\log^p(\mu_i)) \leq \eta^p$, we assume that this happens *in expectation* so that $F$ obeys $\mathcal{E}_F \log(\mu_1)^p \leq \eta^p$. We call this viewpoint the *Bayesian model* because now the estimands are random. Although it seems a digression from our original purposes, it is interesting in its own right and will be connected back to the original model in Section 7.

The motivation for this model is, of course, the ease of analysis. We immediately obtain the asymptotic consistency of FDR thresholding as given in the following:

COROLLARY 4.1. *For $G \in \mathcal{G}$ and $G \neq E$, the empirical FDR threshold $T_q(G_n)$ converges to $T_q(G)$, that is,*

$$\lim_{n \to \infty} T_q(G_n) = T_q(G), \qquad a.s.$$

In a natural sense, the FDR functional $T_q(G)$ can be considered as the *ideal FDR* threshold—the threshold that FDR is "trying" to estimate and use.

PROOF. The 'Fundamental Theorem of Statistics' (for example, [16], page 1) tells us that if $G_n$ is the empirical cdf of $X_1, X_2, \ldots, X_n$ i.i.d. $G$, then

(4.11) $$\|G_n - G\| \to 0, \qquad \text{a.s.}$$

Simply combining this with continuity of $T_q(G)$ at $G \neq E$ gives the proof. $\square$

Of course, we can sharpen our conclusions to rates. Under i.i.d. sampling $X_i \sim G$, $\|G_n - G\| = O_P(n^{-1/2})$. Matching this, we have a root-$n$ rate of convergence for the FDR functional.

COROLLARY 4.2. *If $G \in \mathcal{G}$ and $G \neq E$, then*

$$|T_q(G_n) - T_q(G)| = O_P(n^{-1/2}),$$

*where the $O_P()$ is locally uniform in $G$.*

PROOF. Indeed,

$$|T_q(G_n) - T_q(G)| \leq \omega^*(\|G_n - G\|; T_q(G)) = \omega^*(O_P(n^{-1/2}); T_q(G)).$$

By (4.4), for small $\varepsilon > 0$, $\omega^*(\varepsilon; T_q(G)) \leq C_G \varepsilon$, where $C_G$ locally bounded when $G \neq E$. Therefore, this last term is locally uniformly $O_P(n^{-1/2})$ at each $G \in \mathcal{G}$ where $G \neq E$. $\square$

We can, of course, go further. By Massart's work on the DKW constant [9, 15], we have

(4.12) $$P\{\|G_n - G\| \geq s/\sqrt{n}\} \leq 2e^{-2s^2}, \qquad \forall s \geq 0,$$

which combines with estimates of $\omega^*$ to control probabilities of deviations $T_q(G_n) - T_q(G)$.

**5. Ideal FDR thresholding.** Continuing now in the Bayesian model just defined, we define the *ideal FDR thresholding* pseudo-estimate $\tilde{\mu}_{q,n}$, with coordinates $(\tilde{\mu}_i)$ given by

$$\tilde{\mu}_i = \begin{cases} X_i, & X_i \geq T_q(G), \\ 1, & \text{otherwise.} \end{cases} \tag{5.1}$$

In words, we are thresholding at the large-sample limit of the FDR procedure.

Note that $T_q(G)$ depends on the underlying cdf $G$, which is actually unknown in any realistic situation. $\tilde{\mu}_{q,n}$ is not a *true* estimator; it could only be applied in a setting where we had side information supplied by an *oracle* which told us $T_q(G)$. We view $\tilde{\mu}_{q,n}$ as an *ideal procedure* and the risk for $\tilde{\mu}_{q,n}$ as an *ideal risk*—the risk we would achieve if we could use the threshold that FDR is 'trying' to 'estimate.' Despite the gap between 'true' and 'ideal,' $\tilde{\mu}_{q,n}$ plays an important role in studying the true risk for (true) FDR thresholding. In fact, we will eventually show that, asymptotically, there is only a negligible difference between the ideal risk for $\tilde{\mu}_{q,n}$ and the (true) risk for the FDR thresholding estimator $\hat{\mu}_{q,n}$. Let $\tilde{\mathcal{R}}_n(T_q, G)$ denote the ideal risk for $\tilde{\mu}_{q,n}$ in the Bayesian model,

$$\tilde{\mathcal{R}}_n(T_q, G) \equiv \frac{1}{n} \mathcal{E}\left[\sum_{i=1}^n (\log(\tilde{\mu}_{q,n})_i - \log \mu_i)^2\right].$$

Arguing much as in Sections 2 and 3 above, in the Bayesian model, we also have the following identity with univariate thresholding risk:

$$\tilde{\mathcal{R}}_n(T_q, G) = \rho_T(T_q(G), F). \tag{5.2}$$

Since this ideal risk depends only on a univariate random variable $X_1 \sim G$ and $T_q(G)$ is nonstochastic, its analysis is relatively straightforward. Also, we can now drop the subscript $n$ from $\tilde{\mathcal{R}}_n$.

THEOREM 5.1. *Fix $0 < q < 1$ and $0 < p < 2$.*

1. *Worst-case ideal risk. We have*

$$\lim_{\eta \to 0} \left[\frac{\sup_{G \in \mathcal{G}_p(\eta)} \tilde{\mathcal{R}}(T_q, G)}{\eta^p \log^{2-p} \log \frac{1}{\eta}}\right] = \begin{cases} 1, & 0 < q \leq \frac{1}{2}, \\ \dfrac{q}{1-q}, & \frac{1}{2} < q < 1. \end{cases} \tag{5.3}$$

2. *Least favorable scale mixture. Fix $0 \leq s \leq 1$. Set*

$$\mu_b^* = \mu_b^*(\eta) = \log\left(\frac{1}{\eta}\right) \bigg/ \log\log\left(\frac{1}{\eta}\right), \qquad \mu_v^* = \mu_v^*(\eta) = \log\left(\frac{1}{\eta}\right) \cdot \log\log\left(\frac{1}{\eta}\right)$$

*and*

$$G_{\varepsilon,\mu} = (1-\varepsilon)E(\cdot) + \varepsilon E(\cdot/\mu), \qquad \varepsilon \cdot \log^p(\mu) = \eta^p.$$

*Define*

$$\tilde{\mu} = \tilde{\mu}(\eta; q, s) = \begin{cases} \mu_b^*(\eta), & 0 < q < \frac{1}{2}, \\ \mu_v^*(\eta), & \frac{1}{2} < q < 1, \\ (1-s) \cdot \mu_b^*(\eta) + s \cdot \mu_v^*(\eta), & 0 \le s \le 1, \quad q = \frac{1}{2}. \end{cases}$$

*Then $G_{\varepsilon,\tilde{\mu}}$ is asymptotically least favorable for $T_q$, that is,*

$$\lim_{\eta \to 0} \left[ \frac{\tilde{\mathcal{R}}(T_q, G_{\varepsilon,\tilde{\mu}})}{\sup_{G \in \mathcal{G}_p(\eta)} \tilde{\mathcal{R}}(T_q, G)} \right] = 1.$$

By Theorems 2.1–2.3, the denominator on the left-hand side of (5.3) is asymptotically equivalent to the minimax risk in the original model of Section 1. In words, the worst-case ideal risk for the i.i.d. sampling model is asymptotically equivalent to the minimax risk (1.4) as $\eta \to 0$. This, of course, is no accident; it is a key step towards Theorem 1.3.

5.1. *Proof of Theorem 5.1.* We now describe, in a series of lemmas the ideas for proving Theorem 5.1. In later subsections, we prove the individual lemmas.

Since the ideal risk $\tilde{\mathcal{R}}(T_q, G)$ is, by (5.2), reducible to the univariate thresholding Bayes risk which we studied in Section 3, we know to split the ideal risk $\tilde{\mathcal{R}}(T_q, G)$ into two terms, the *bias* proxy and the *variance* proxy,

$$\tilde{B}^2(T_q, G) \equiv \int b(T_q(G), \mu) \, dF(\mu), \qquad \tilde{V}(T_q, G) \equiv \int v(T_q(G), \mu) \, dF(\mu).$$

Consider $\tilde{V}(T_q, G)$. Asymptotically as $\eta \to 0$, every eligible $F \in \mathcal{F}_p(\eta)$ puts almost all mass in the vicinity of 1, so

(5.4) $$\tilde{V}(T_q, G) \approx v(T_q(G), 1) \approx \log^2(T_q(G)) e^{-T_q(G)}.$$

We set $\tilde{v}(t) \equiv \log^2(t) e^{-t}$. The following formal approximation result is proved in [11], Chapter 6:

LEMMA 5.1. *As $\eta \to 0$,*

$$\sup_{G \in \mathcal{G}_p(\eta)} |\tilde{V}(T_q, G) - \tilde{v}(T_q(G))| = o\left( \eta^p \log^{2-p} \log \frac{1}{\eta} \right).$$

Note that as $G$ tends to $E$, Lemma 4.4 implies that $T_q(G) \to \infty$. Since $\tilde{v}(T_q(G))$ decreases rapidly, the key to majorizing the variance is to keep $T_q(G)$ small, motivating study of

(5.5) $$T_q^* = T_q^*(\eta; p) = \inf_{G \in \mathcal{G}_p(\eta)} T_q(G).$$

LEMMA 5.2. *As $\eta \to 0$,*
$$T_q^* = T_q^*(\eta; p) = p\left(\log\frac{1}{\eta} + \log\log\log\frac{1}{\eta}\right) + \log\left(\frac{1-q}{q}\right) + o(1).$$

The proof is given in Section 5.2. As a direct result, we get
$$\log^2(T_q^*)e^{-T_q^*} = \left[\frac{q}{1-q}\eta^p \log^{2-p}\frac{1}{\eta}\right] \cdot (1 + o(1));$$

moreover, when $T_q(G)$ exceeds $T_q^*$, the variance proxy $\tilde{v}(T_q^*)$ drops and we obtain the following:

LEMMA 5.3. *As $\eta \to 0$,*
$$\sup_{G \in \mathcal{G}_p(\eta)} \tilde{V}(T_q, G) = \left[\frac{q}{1-q}\eta^p \log^{2-p}\log\frac{1}{\eta}\right] \cdot (1 + o(1))$$

*and*
$$\sup_{G \in \mathcal{G}_p(\eta), T_q(G) \geq T_q^* + \sqrt{T_q^*}} \tilde{V}(T_q, G) = o\left(\eta^p \log^{2-p}\log\frac{1}{\eta}\right).$$

We now study the bias proxy. The key observation is as follows:

(5.6) $$b(t, \mu) \approx \begin{cases} \log^2 \mu, & \mu \ll t, \\ \frac{t}{\mu}\log^2 \mu, & \mu \gg t. \end{cases}$$

To develop intuition, consider the family of 2-point mixtures
$$\mathcal{G}_p^{2,0}(\eta) = \{G_{\varepsilon,\mu} = (1-\varepsilon)E(\cdot) + \varepsilon E(\cdot/\mu), \varepsilon \log^p \mu = \eta^p\}.$$

Now, (5.6) tells us that the maximum of the bias functional over this family is obtained by taking $\mu$ as large as possible, while avoiding $\frac{T_q(G_{\varepsilon,\mu})}{\mu} \ll 1$; moreover, direct calculations show that

(5.7) $$\frac{T_q(G_{\varepsilon,\mu})}{\mu} = \frac{\log(1 + p(\frac{1}{q}-1)\frac{1}{\eta^p}\log(\mu))}{\mu - 1},$$

so the value of $\mu$ causing the worst bias proxy should be close to the solution of the following equation:
$$\frac{\log(1 + p(\frac{1}{q}-1)\frac{1}{\eta^p}\log(\mu))}{\mu - 1} = 1.$$

Elaborating on this idea leads to the following result, to be proven in Section 5.3:

LEMMA 5.4. *As $\eta \to 0$,*

$$\sup_{G \in \mathcal{G}_p(\eta)} \tilde{B}^2(T_q, G) = \left( \eta^p \log^{2-p} \log \frac{1}{\eta} \right) \cdot (1 + o(1)).$$

Combine the above analysis for bias and variance proxies, to give

$$1 + o(1) \leq \frac{\sup_{G \in \mathcal{G}_p(\eta)} \tilde{\mathcal{R}}(T_q, G)}{\eta^p \log^{2-p} \log \frac{1}{\eta}} \leq \frac{1}{1-q} + o(1) \qquad \text{as } \eta \to 0.$$

Compare this to the conclusion of Theorem 5.1; we have obtained the correct rate, but not yet the precise constant. To refine our analysis, note that the worst bias and the worst variance are obtained at different values $\mu$ within the family $\mathcal{G}_p^{2,0}(\eta)$. Denote the $\mu$'s causing the worst bias and the worst variance by $\mu_b^*$ and $\mu_v^*$. Then

$$\mu_b^* \sim \frac{\log \frac{1}{\eta}}{\log \log \frac{1}{\eta}}, \quad \mu_v^* \sim \log \frac{1}{\eta} \cdot \log \log \frac{1}{\eta} \qquad \text{as } \eta \to 0.$$

Divide $\mathcal{G}_p(\eta)$ into two subsets,

$$\mathcal{G}_1 \equiv \{ G \in \mathcal{G}_p(\eta), T_q(G) \geq T_q^* + \sqrt{T_q^*} \},$$

$$\mathcal{G}_2 \equiv \{ G \in \mathcal{G}_p(\eta), T_q(G) < T_q^* + \sqrt{T_q^*} \}$$

and consider each separately. [Note that $G_{\mu_b^*} \in \mathcal{G}_1$, while $G_{\mu_v^*} \in \mathcal{G}_2$. Here, $G_{\mu_b^*}$ and $G_{\mu_v^*}$ are mixtures of point masses at 1 and $\mu$ living in $\mathcal{G}_p^{2,0}(\eta)$ with $\mu = \mu_b^*$ and $\mu_v^*$, respectively]. Over the first subset, the variance is uniformly $O(\eta^p)$ and we immediately obtain

$$\sup_{\mathcal{G}_1} \tilde{\mathcal{R}}(T_q, G) \approx \sup_{\mathcal{G}_1} \tilde{B}^2(T_q, G) \approx \eta^p \log^{2-p} \log \frac{1}{\eta} \qquad \text{as } \eta \to 0.$$

For the second subset, the following lemma is proved in [6], page 22:

LEMMA 5.5. *As $\eta \to 0$,*

$$\sup_{\mathcal{G}_2} \tilde{\mathcal{R}}(T_q, G) = \begin{cases} \left( \eta^p \log^{2-p} \log \frac{1}{\eta} \right) \cdot (1 + o(1)), & 0 < q \leq \frac{1}{2}, \\ \frac{q}{1-q} \cdot \left( \eta^p \log^{2-p} \log \frac{1}{\eta} \right) \cdot (1 + o(1)), & \frac{1}{2} < q < 1. \end{cases}$$

Theorem 5.1 follows once Lemmas 5.2 and 5.4 are proved.

5.2. *Proof of Lemma 5.2.* Consider the upper envelope of the survivor function among moment-constrained scale mixtures,

$$\bar{G}_t^* = \bar{G}_t^*(\eta; p) = \sup\{\bar{G}(t), G \in \mathcal{G}_p(\eta)\}.$$

The quantity of interest is the crossing point where this envelope meets the FDR boundary,

$$T_q^* = \inf\{t : \bar{G}_t^* \geq \bar{E}(t)/q\}.$$

Equivalently,

(5.8) $$T_q^* = \inf\{t : [(\bar{G}_t^*/\bar{E}(t)) - 1] \geq (1-q)/q\}.$$

Letting

$$h^*(t; \eta, p) = [(\bar{G}_t^*/\bar{E}(t)) - 1],$$

the key to calculating $T_q^*$ is to explicitly express $h^*(t)$ as a function of $t$, asymptotically, for small $\eta$.

Calculating $h^*(t)$ again involves optimization of a linear functional over a class of moment-constrained cdf's and we can apply the theory in Section 3.2. Set $\psi = \psi_t(\mu) = [e^{(1-\frac{1}{\mu})t} - 1]$ and $\phi(\mu) = \log^p(\mu)$ and define $\Psi = \Psi_t$ as in (3.10) so that $h^*(t; \eta, p) = \Psi_t(\eta^p)$. Note that

(5.9) $$\lim_{\mu \to 1+} \left[\frac{\psi_t(\mu)}{\log^p(\mu)}\right] = \begin{cases} 0, & 0 < p < 1, \\ t, & p = 1, \\ \infty, & 1 < p < 2, \end{cases}$$

so we treat the cases $0 < p \leq 1$ and $1 < p < 2$ separately.

When $0 < p \leq 1$, elementary analysis shows that for large $t$,

$$\mu^* = \underset{\mu \geq 1}{\operatorname{argmax}}\left\{\frac{e^{(1-\frac{1}{\mu})t} - 1}{\log^p(\mu)}\right\} \sim \frac{t \log(t)}{p}, \qquad \Psi^* = \frac{e^{(1-\frac{1}{\mu^*})t} - 1}{\log^p(\mu^*)} \sim e^t/[\log^p(t)],$$

so the condition of Lemma 3.1 is satisfied and

(5.10) $$\Psi_t(\eta^p) \sim \eta^p e^t / \log^p(t).$$

Inserting (5.10) into (5.8) and solving for $t$ gives the lemma for the case $0 < p \leq 1$.

When $1 < p < 2$, direct calculations show that the function $\psi'(\mu)/\phi'(\mu)$ strictly increases in the interval $(1, \bar{\mu}]$ with $\log(\bar{\mu}) = \log(\bar{\mu}(t; p)) = (p-1)/t$, also that $[\psi'(\bar{\mu})/\phi'(\bar{\mu})] \leq \Psi^{**}(\bar{\mu})$, so the condition of Lemma 3.2 is satisfied. More calculations show first, that,

$$\mu^* = \mu^*(t; p) \sim \underset{\{\mu' \geq \bar{\mu}\}}{\operatorname{argmax}}\left\{\frac{\psi(\mu')}{\log^p(\mu')}\right\} \sim \frac{t}{p \log(t)},$$

second, that for any $1 < \mu \leq \bar{\mu}$,

$$\Psi^{**}(\mu) = \Psi^{**}(\mu;t)$$
$$\equiv \max_{\{\mu' \geq \bar{\mu}\}} \left\{ \frac{\psi(\mu') - \psi(\mu)}{\log^p(\mu') - \log^p(\mu)} \right\} \sim \max_{\{\mu' \geq \bar{\mu}\}} \left\{ \frac{\psi(\mu')}{\log^p(\mu')} \right\} \sim \frac{e^t}{\log^p(t)}$$

and, finally, that

$$\log(\mu_*) = \log(\mu_*(t;p)) \sim \left( \frac{1}{p} t \log^p(t) e^{-t} \right)^{1/(p-1)}$$

since $h^*(t, \eta, p) = \Psi_t(\eta^p)$. By Lemma 3.2,

$$(5.11) \quad h^*(t, \eta, p) = \begin{cases} e^{(1-e^{-\eta})t} - 1, & \eta^p \leq \log^p(\mu_*), \\ e^{(1-\frac{1}{\mu_*})t} - 1 + \Psi^{**}(\mu_*)(\eta^p - \log(\mu_*)), \\ & \log^p(\mu_*) < \eta^p \leq \log^p(\mu^*); \end{cases}$$

moreover, by letting $t^* = t_p^*(\eta)$ denote the solution of $\log^p(\mu_*(t,p)) = \eta^p$, we can rewrite (5.11) as

$$(5.12) \quad h^*(t; \eta, p) = \begin{cases} e^{(1-e^{-\eta})t} - 1, & t \leq t^*, \\ e^{(1-\frac{1}{\mu_*})t} - 1 + \Psi^{**}(\mu_*)(\eta^p - \log(\mu_*)), & t \geq t^*, \end{cases}$$

here noting that $t^* \sim (p-1) p \log(\frac{1}{\eta})$ for small $\eta$.

Inserting (5.12) into (5.8), it becomes clear that for sufficiently small $\eta$ and $t \leq t^*$, $h(t; \eta, p) \approx 0$. Thus, $T_q^*$ is obtained by equating

$$\frac{1-q}{q} = e^{(1-\frac{1}{\mu_*})t} - 1 + \Psi^{**}(\mu_*)(\eta^p - \log(\mu_*)) \sim \eta^p e^t / \log^p(t),$$

which gives the lemma for the case $1 < p < 2$. $\square$

5.3. *Proof of Lemma 5.4.*

LEMMA 5.6. *For a measurable function $\psi$ defined on $[1, \infty)$, where $\psi \geq 0$ but is not identically 0 and $\sup_{\mu \geq 1} \{\psi(\mu)/\mu\} < \infty$, then for $G \in \mathcal{G}$ and $0 < \tau < T_q(G)$, we have*

$$\int \psi(\mu) [e^{-\tau/\mu} - e^{-T_q(G)/\mu}] dF \leq (1/q) \sup_{\{\mu \geq 1\}} \{\psi(\mu)/\mu\} \cdot \tau e^{-\tau}/(1 - e^{-\tau}).$$

Letting $\tau \to 0$ and combining Lemma 5.6 with Fatou's Lemma, we have

$$(5.13) \quad \int \psi(\mu) [1 - e^{-T_q(G)/\mu}] dF \leq (1/q) \sup_{\{\mu \geq 1\}} \{\psi(\mu)/\mu\}.$$

PROOF. Let $k_0 = k_0(\tau; G) = \lfloor \frac{T_q(G)}{\tau} \rfloor$. Since $T_q(G) > \tau$, $k_0 \geq 1$. Moreover,

$$(5.14) \quad \int \psi(\mu)[e^{-\tau/\mu} - e^{-T_q(G)/\mu}] \, dF \leq \int \psi(\mu)[e^{-\tau/\mu} - e^{-(k_0+1)\tau/\mu}] \, dF$$

$$(5.15) \quad = \int \psi(\mu)(1 - e^{-\tau/\mu}) \left[ \sum_{j=1}^{k_0} e^{-j \cdot \tau/\mu} \right] dF.$$

We introduce the shorthand notation $c = \max_{\mu \geq 1}\{\psi(\mu)/\mu\}$ and recall that $1 - e^{-x/\mu} \leq x/\mu$ for all $x \geq 0$, so for $1 \leq j \leq k_0$,

$$(5.16) \quad \int \psi(\mu)(1 - e^{-\tau/\mu}) e^{-j \cdot \tau/\mu} \, dF \leq \tau \int (\psi(\mu)/\mu) e^{-j \cdot \tau/\mu} \, dF$$

$$\leq \tau \cdot c \cdot \int e^{-j \cdot \tau/\mu} \, dF.$$

By definition of $k_0$ and the FDR functional,

$$(5.17) \quad \int e^{-j \cdot \tau/\mu} \, dF = \bar{G}(j \cdot \tau) \leq (1/q) e^{-j \cdot \tau}, \qquad 1 \leq j \leq k_0.$$

Combining (5.14)–(5.17) gives

$$(5.18) \quad \int \psi(\mu)[e^{-\tau/\mu} - e^{-T_q(G)/\mu}] \, dF \leq (c/q) \cdot \tau \cdot \sum_{j=1}^{k_0} e^{-j \cdot \tau}$$

$$\leq (c/q) \cdot \tau \cdot e^{-\tau}/(1 - e^{-\tau}). \qquad \square$$

We now prove Lemma 5.4. As in Section 3, let

$$t_0 = t_0(p, \eta) = p \log(1/\eta) + p \log\log(1/\eta) + \sqrt{\log\log(1/\eta)}.$$

By the monotonicity of $b(t, \mu)$ and (3.8), for sufficiently small $\eta > 0$,

$$(5.19) \quad \sup_{G \in \mathcal{G}_p(\eta), T_q(G) \leq t_0} \tilde{B}^2(T_q, G) \leq \sup_{\mathcal{G}_p(\eta)} \int b(t_0, \mu) \, dF$$

$$= \eta^p \log^{2-p} \log(1/\eta)(1 + o(1)).$$

Moreover, for any $G$ with $T_q(G) > t_0$, letting $\psi(\cdot) = \log^2(\cdot)$ and $\tau = t_0$ in Lemma 5.6, we have

$$0 \leq \tilde{B}^2(T_q, G) - \int b(t_0, \mu) \, dF = \int \log^2(\mu)[e^{-t_0/\mu} - e^{-T_q(G)/\mu}] \, dF$$

$$\leq c t_0 e^{-t_0}/(1 - e^{-t_0}),$$

where $c = \max_{\mu \geq 1}\{\log^2(\mu)/\mu\}$, so it is clear that

$$(5.20) \quad \sup_{\{G \in \mathcal{G}_p(\eta), T_q(G) > t_0\}} \tilde{B}^2(T_q, G) \leq \int b(t_0, \mu) \, dF + O(t_0 e^{-t_0}).$$

Lemma 5.4 follows directly from (5.19)–(5.20) and $t_0 e^{-t_0} = o(\eta^p \log^{2-p} \log(\frac{1}{\eta}))$.
$\square$

**6. Asymptotic risk behavior for FDR thresholding.** We now turn to $\hat{\mu}_{q,n}$, the true FDR thresholding estimator. For technical reasons, we define a threshold $\hat{T}_{q,n}$ slightly differently than $\hat{t}^{FDR}$. This difference does not affect the estimate. Thus, we will have $\hat{\mu}_{q,n} \equiv \hat{\mu}_{\hat{T}_{q,n}} = (\hat{\mu}_i)$ with

$$\hat{\mu}_i = \begin{cases} X_i, & X_i \geq \hat{T}_{q,n}, \\ 1, & X_i < \hat{T}_{q,n}. \end{cases}$$

Our strategy is to show that the ideal and true FDR behave similarly.

Still in the Bayesian model, we let $\mathcal{R}_n(\hat{T}_{q,n}, G)$ denote the per-coordinate average risk for $\hat{\mu}_{q,n}$, that is,

$$\mathcal{R}_n(\hat{T}_{q,n}, G) \equiv \frac{1}{n}\mathcal{E}\left[\sum_{i=1}^n (\log(\hat{\mu}_{q,n})_i - \log\mu_i)^2\right].$$

Here, again, the expectation is over $(X_i, \mu_i)$ pairs i.i.d. with bivariate structure $X_i|\mu_i \sim \text{Exp}(\mu_i)$.

We will show that as $n \to \infty$, the difference between the true risk $\mathcal{R}_n(\hat{T}_{q,n}, G)$ and the ideal risk $\tilde{\mathcal{R}}(T_q, G)$ is asymptotically negligible. We suppress the subscript $n$ on $\mathcal{R}_n$ (this is an abuse of notation).

THEOREM 6.1.

$$\lim_{n\to\infty}\left[\sup_{G\in\mathcal{G}}|\mathcal{R}(\hat{T}_{q,n}, G) - \tilde{\mathcal{R}}(T_q, G)|\right] = 0.$$

As a result,

$$\lim_{n\to\infty}\left[\sup_{G\in\mathcal{G}_p(\eta)}|\mathcal{R}(\hat{T}_{q,n}, G) - \tilde{\mathcal{R}}(T_q, G)|\right] = 0.$$

Combining Theorems 6.1 and 5.1, we have

$$\lim_{\eta\to 0}\left[\lim_{n\to\infty}\frac{\sup_{G\in\mathcal{G}_p(\eta)}\mathcal{R}(\hat{T}_{q,n}, G)}{\eta^p \log^{2-p}\log\frac{1}{\eta}}\right] = \begin{cases} 1, & 0 < q \leq \frac{1}{2}, \\ \dfrac{q}{1-q}, & \frac{1}{2} < q < 1. \end{cases}$$

Hence, $\hat{T}_{q,n}$ asymptotically achieves the $n$-variate minimax Bayes risk when $n \to \infty$ followed by $\eta \to 0$.

6.1. *Proof of Theorem 6.1.* We begin by defining $\hat{T}_{q,n}$. In applying the FDR functional to the empirical distribution, it is always possible that

(6.1) $$\bar{G}_n(t) < \frac{1}{q}\bar{E}(t), \quad \text{for all } t > 0,$$

in which case $T_q(G_n) = \hat{t}^{FDR} = +\infty$. Letting $W_n$ denote the event (6.1), define

$$\hat{T}_{q,n} = \begin{cases} T_q(G_n), & \text{over } W_n^c, \\ \log(\frac{n}{q}), & \text{over } W_n. \end{cases} \quad (6.2)$$

The following lemma, which was proven in [6, 11], shows that this definition of threshold gives the same estimator as $T_q(G_n)$, while obeying a bound which is convenient for analysis:

LEMMA 6.1. *Suppose $X_i \overset{iid}{\sim} G$, $G \in \mathcal{G}$, $G \neq E$ and $\hat{T}_{q,n}$ is defined as in (4.1). Then*

1. *The FDR estimator is equivalently realized by thresholding at $\hat{T}_{q,n}$:*
$\hat{\mu}_{q,n}^{FDR} = \hat{\mu}_{\hat{T}_{q,n}}$.
2. $\hat{T}_{q,n} \leq \log(\frac{n}{q})$.

Next, we study the risk for $\hat{T}_{q,n}$. We have

$$\mathcal{R}(\hat{T}_{q,n}, G) = \frac{1}{n}\sum_{i=1}^{n} \mathcal{E}_F \mathcal{E}_\mu \left[ \log^2(\mu_i) \mathbb{1}_{\{X_i < \hat{T}_{q,n}\}} + \log^2\left(\frac{X_i}{\mu_i}\right) \mathbb{1}_{\{X_i \geq \hat{T}_{q,n}\}} \right]$$

$$= \mathcal{E}_F \mathcal{E}_\mu [\log^2(\mu_1) \mathbb{1}_{\{X_1 < \hat{T}_{q,n}\}} + \log^2(X_1/\mu_1) \mathbb{1}_{\{X_1 \geq \hat{T}_{q,n}\}}]$$

and $\mathcal{R}(\hat{T}_{q,n}, G)$ naturally splits into a 'bias' proxy and the 'variance' proxy, as follows:

$$B^2(\hat{T}_{q,n}, G) = \mathcal{E}_F \mathcal{E}_\mu [\log^2(\mu_1) \mathbb{1}_{\{X_1 < \hat{T}_{q,n}\}}],$$

$$V(\hat{T}_{q,n}, G) = \mathcal{E}_F \mathcal{E}_\mu [\log^2(X_1/\mu_1) \mathbb{1}_{\{X_1 \geq \hat{T}_{q,n}\}}].$$

The comparable notions in the ideal risk case were

$$\tilde{B}^2(T_q, G) = \mathcal{E}_F \mathcal{E}_\mu [\log^2(\mu_1) \mathbb{1}_{\{X_1 < T_q(G)\}}],$$

$$\tilde{V}(T_q, G) = \mathcal{E}_F \mathcal{E}_\mu [\log^2(X_1/\mu_1) \mathbb{1}_{\{X_1 \geq T_q(G)\}}].$$

Intuitively, we expect that $\tilde{B}^2$ is 'close' to $B^2$ and that $\tilde{V}$ is 'close' to $V$; our next task is to validate these expectations. Observe that

$$(6.3) \quad |B^2(\hat{T}_{q,n}, G) - \tilde{B}^2(T_q, G)| \leq \mathcal{E}[\log^2(\mu_1) |\mathbb{1}_{\{X_1 < \hat{T}_{q,n}\}} - \mathbb{1}_{\{X_1 < T_q(G)\}}|],$$

$$(6.4) \quad |V(\hat{T}_{q,n}, G) - \tilde{V}(T_q, G)| \leq \mathcal{E}[\log^2(X_1/\mu_1) |\mathbb{1}_{\{X_1 < \hat{T}_{q,n}\}} - \mathbb{1}_{\{X_1 < T_q(G)\}}|].$$

It would not be hard to validate the expectations if $|\hat{T}_{q,n} - T_q(G)|$ were negligible for large $n$, uniformly for $G \in \mathcal{G}$. In Section 4, Lemma 4.5 tells us that $\hat{T}_{q,n} - T_q(G)$ is locally $O_P(n^{-1/2})$ or, more specifically,

(6.5) $\quad |T_q(G) - T_q(G_n)| \sim \dfrac{q}{\log(1/q)} T_q(G) e^{T_q(G)} \|G - G_n\|, \qquad G \neq E.$

Unfortunately, for any fixed $n$, $G$ might get arbitrary close to $E$ and, as a result, $T_q(G)$ might get arbitrary large, so the relationship in (6.5) cannot hold *uniformly* over $G \in \mathcal{G}$.

A closer look reveals that those $G$'s failing (6.5) would, roughly, satisfy

$$T_q(G) e^{T_q(G)} \geq \sqrt{n}, \quad \text{or} \quad T_q(G) \geq \log(n)/2.$$

Note that as $n$ increases from 1 to $\infty$, $\{G \in \mathcal{G} : T_q(G) \geq \log(n)/2\}$ defines a sequence of subsets, strictly decreasing to $\varnothing$. Motivated by this, we look for a subsequence of subsets of $\mathcal{G}$ obeying

(a) $\mathcal{G}^{(1)} \subset \mathcal{G}^{(2)} \subset \cdots \subset \mathcal{G}^{(n)} \subset \cdots$ and $\bigcup_1^\infty \mathcal{G}^{(n)} = \mathcal{G}$;
(b) $\mathcal{G}^{(n)}$ approaches $\mathcal{G}$ *slowly* enough such that $\sup_{\mathcal{G}^{(n)}} [\sqrt{n} T_q(G) e^{T_q(G)}] = o(1)$;
(c) for large $n$, $|\mathcal{R}(\hat{T}_{q,n}) - \tilde{\mathcal{R}}(T_q, G)|$ is uniformly negligible over $\mathcal{G} \setminus \mathcal{G}^{(n)}$.

A convenient choice is

(6.6) $\qquad \mathcal{G}_1^{(n)} \equiv \{G \in \mathcal{G} : T_q(G) \leq \log(n)/8\}, \qquad n \geq 1.$

We expect that the difference between $T_q(G_n)$ and $T_q(G)$ is uniformly negligible over $\mathcal{G}_1^{(n)}$, that is,

$$\sup_{\mathcal{G}_1^{(n)}} |T_q(G) - T_q(G_n)| = o_p(1).$$

LEMMA 6.2. *Let $A_n$ denote the event $\{|\hat{T}_{q,n} - T_q(G)| \leq n^{-1/4}\}$. Then for sufficiently large $n$,*

$$\sup_{G \in \mathcal{G}_1^{(n)}} P_G\{A_n^c\} \leq 3 e^{-[32(1-q)^2/q^2] n^{1/4}/\log^2(n)}.$$

Based on Lemma 6.2, one can develop a proof for the following:

LEMMA 6.3. *For sufficiently small $0 < \delta < 1$,*

1. $\lim_{n \to \infty} \sup_{G \in \mathcal{G}_1^{(n)}} |B^2(\hat{T}_{q,n}, G) - \tilde{B}^2(T_q, G)| = 0$;
2. $\lim_{n \to \infty} \sup_{G \in \mathcal{G}_1^{(n)}} |V(\hat{T}_{q,n}, G) - \tilde{V}(T_q, G)| = 0.$

*As a result, $\lim_{n \to \infty} \sup_{G \in \mathcal{G}_1^{(n)}} |\mathcal{R}(\hat{T}_{q,n}, G) - \tilde{\mathcal{R}}(T_q, G)| = 0.$*

We now consider (c). Define

(6.7) $$\mathcal{G}_0^{(n)} \equiv \mathcal{G} \setminus \mathcal{G}_1^{(n)}, \qquad n \geq 1.$$

Though it is no longer sensible to require that $|T_q(G_n) - T_q(G)|$ be uniformly negligible over $\mathcal{G}_0^{(n)}$, we still hope that $T_q(G_n)$ at least stays at the *same* magnitude as $T_q(G)$, or $T_q(G_n) = O_p(\log(n))$. This turns out to be true and, in fact, is an immediate consequence of Massart's inequality (4.12).

LEMMA 6.4. *Letting $D_n$ be the event $\{\hat{T}_{q,n} \geq \log(n)/16\}$,*

$$\sup_{G \in \mathcal{G}_0^{(n)}} P_G\{D_n^c\} = 2e^{-2[(1-\sqrt{q})^2/q^2]n^{7/8}}.$$

Combining this with Lemma 6.1, we have, except for an event with negligible probability,

$$\log(n)/16 \leq \hat{T}_{q,n} \leq \log(n/q).$$

Since $v(t,\mu)$ is monotone decreasing in $t$, it is now clear that both $V(\hat{T}_{q,n}, G)$ and $\tilde{V}(T_q, G)$ are uniformly negligible over $\mathcal{G}_0^{(n)}$.

LEMMA 6.5.

$$\lim_{n\to\infty}\left[\sup_{G\in\mathcal{G}_0^{(n)}} \tilde{V}(T_q, G)\right] = 0, \qquad \lim_{n\to\infty}\left[\sup_{G\in\mathcal{G}_0^{(n)}} V(\hat{T}_{q,n}, G)\right] = 0.$$

Finally, note that $b(t,\mu)$ is strictly increasing in $t$, so either $B^2(\hat{T}_{q,n}, G)$ or $\tilde{B}^2(T_q, G)$ will not be uniformly negligible over $\mathcal{G}_0^{(n)}$. However, note that $b(t,\mu)$ increases very *slowly* in $t$ for large $t$, so we can expect that $|B^2(\hat{T}_{q,n}, G) - \tilde{B}^2(T_q, G)|$ is uniformly negligible over $\mathcal{G}_0^{(n)}$.

LEMMA 6.6. $\lim_{n\to\infty}[\sup_{G\in\mathcal{G}_0^{(n)}} |B^2(\hat{T}_{q,n}, G) - \tilde{B}^2(T_q, G)|] = 0.$

The choice of $\log(n)/8$ is only for convenience; a similar result holds if we replace $\log(n)/8$ by $c\log(n)$ for $0 < c < 1/2$.

Combining the above lemmas yields Theorem 6.1. □

The proofs of Lemmas 6.1–6.6 can be found in the full version of this paper [6].

**7. Proof of Theorem 1.3.** We now complete the proof of Theorem 1.3. The key point is to relate the Bayesian model of Sections 4–6 to the frequentist model of Section 1. In the frequentist model, $X_i \sim \text{Exp}(\mu_i), 1 \leq i \leq n$, where $\mu = (\mu_1, \mu_2, \ldots, \mu_n)$ is an arbitrary deterministic vector $\mu \in M_{n,p}(\eta)$. Recall that $\mathcal{R}_n(\hat{T}_{q,n}, G)$ denotes the risk of FDR estimation in the Bayesian model, while $R_n(\hat{\mu}_{q,n}, \mu)$ denotes the risk in the frequentist model. Below, we will show that

$$(7.1) \qquad \lim_{\eta \to 0}\left[\lim_{n\to\infty} \frac{\sup_{G\in\mathcal{G}_p(\eta)} \mathcal{R}_n(\hat{T}_{q,n}, G)}{\sup_{\mu \in M_{n,p}(\eta)} R_n(\hat{\mu}_{q,n}, \mu)}\right] = 1.$$

Recall that by Theorems 1.1, 5.1 and 6.1, we have

$$\lim_{\eta\to 0}\left[\lim_{n\to\infty}\frac{\sup_{G\in\mathcal{G}_p(\eta)}\mathcal{R}_n(\hat{T}_{q,n},G)}{R_n^*(M_{n,p}(\eta))}\right] = \begin{cases} 1, & 0 < q \leq \frac{1}{2}, \\ \dfrac{q}{1-q}, & \frac{1}{2} < q < 1, \end{cases}$$

so Theorem 1.3 follows from (7.1). To prove (7.1), let $G_\mu$ denote the mixture $G_\mu = \frac{1}{n}\sum_{i=1}^{n} E(\cdot/\mu_i)$, let $\tilde{R}_n(\tilde{\mu}_{q,n}, \mu)$ denote the ideal risk for thresholding at $T_q(G_\mu)$ under the frequentist model and let $\tilde{\mathcal{R}}(T_q, G)$ again denote the ideal risk for thresholding at $T_q(G)$ in the Bayesian model. We have the following crucial identity:

$$(7.2) \qquad \tilde{R}_n(\tilde{\mu}_{q,n}, \mu) \equiv \tilde{\mathcal{R}}(T_q, G_\mu), \qquad \forall \mu, n.$$

Also, note that the class of $G_\mu$'s arising from some $\mu \in M_{n,p}(\eta)$ is a subset of the class of all $G$'s arising in $\mathcal{G}_p(\eta)$, for each $n > 0$. Hence,

$$\sup_{\mu \in M_{n,p}(\eta)} \tilde{\mathcal{R}}(T_q, G_\mu) \leq \sup_{G \in \mathcal{G}_p(\eta)} \tilde{\mathcal{R}}(T_q, G), \qquad \forall n.$$

However, note that by Theorem 5.1, appropriately chosen 2-point priors can be asymptotically least-favorable for ideal risk in the Bayesian model. By choosing $\mu$ which contain entries with only the two underlying values in the least favorable prior and with appropriate underlying frequencies, we can obtain

$$(7.3) \qquad \lim_{\eta\to 0}\left[\frac{\lim_{n\to\infty}\sup_{\mu \in M_{n,p}(\eta)} \tilde{\mathcal{R}}(T_q, G_\mu)}{\sup_{G\in\mathcal{G}_p(\eta)}\tilde{\mathcal{R}}(T_q, G)}\right] = 1.$$

Now, relating the Bayesian to the frequentist model via (7.2), we have

$$(7.4) \qquad \lim_{\eta\to 0}\left[\frac{\lim_{n\to\infty}\sup_{\mu\in M_{n,p}(\eta)} \tilde{R}_n(\hat{\mu}_{q,n}, \mu)}{\sup_{G\in\mathcal{G}_p(\eta)}\tilde{\mathcal{R}}(T_q, G)}\right] = 1.$$

Suppose we can next show that the ideal FDR risk in the frequentist model is equivalent to the true risk in the frequentist model, in the same sense as was proved in Theorem 6.1. Hence,

$$\lim_{\eta \to 0} \lim_{n \to \infty} \left[ \frac{\sup_{\mu \in M_{n,p}(\eta)} R_n(\hat{\mu}_{q,n}, \mu)}{\sup_{\mu \in M_{n,p}(\eta)} \tilde{R}_n(\tilde{\mu}_{q,n}, \mu)} \right] = 1. \tag{7.5}$$

Then (7.3)–(7.5) yield (7.1).

The key point is that (7.5) follows exactly as in Section 6. Indeed, there is a precise analog of Theorem 6.1 for the relation between the frequentist risk and the frequentist ideal risk. This is based on two ideas.

First, if $G_n$ now denotes the cdf of $X_1, \ldots, X_n$ *in the frequentist model*, we again have very strong convergence properties of $G_n$, this time to $G_\mu$. This concerns convergence of the empirical cdf for *non-i.i.d.* samples, which is not well known, but can be found in [16], Chapter 25.

LEMMA 7.1 (Bretagnolle [5]). *Let $X_{n1}, X_{n2}, \ldots, X_{nn}$ be independent random variables with arbitrary df's $F_{ni}$, let $F_n(x)$ be the empirical cdf and let $\bar{F} = \text{Ave}_i \{F_{ni}\}$. Then for all $n \geq 1$, $s > 0$, there exists an absolute constant $c$ such that*

$$Prob\{\sqrt{n}\|F_n - \bar{F}_n\| \geq s\} \leq 2ece^{-2s^2}.$$

By means of Massart's work ([16], Chapter 25 and [15]), we can take $c = 1$. Then taking $F_{ni} = \text{Exp}(\mu_i)$ and $\bar{F} = G_\mu$, we obtain

$$P_\mu\{\|G_n - G_\mu\| \geq s/\sqrt{n}\} \leq 6e^{-2s^2}, \qquad \forall \mu.$$

This is completely parallel to the bound (4.12).

Second, it follows immediately from Section 4's analysis that there are frequentist fluctuation bounds for $T_q(G_n) - T_q(G_\mu)$ paralleling those in the Bayesian case. To apply this, we define

$$M_{n,p}^1(\eta) = \{\mu \in M_{n,p}(\eta), T_q(G_\mu) \leq \log(n)/8\} \tag{7.6}$$

and

$$M_{n,p}^0(\eta) = M_{n,p}(\eta) \setminus M_{n,p}^1(\eta). \tag{7.7}$$

LEMMA 7.2. *For sufficiently small $\eta > 0$,*

1. $\lim_{n \to \infty} [\sup_{\mu \in M_{n,p}^1(\eta)} |R_n(\hat{\mu}_{q,n}, \mu) - \tilde{R}_n(\tilde{\mu}_{q,n}, \mu)|] = 0;$
2. $\lim_{n \to \infty} [\sup_{\mu \in M_{n,p}^0(\eta)} |R_n(\hat{\mu}_{q,n}, \mu) - \tilde{R}_n(\tilde{\mu}_{q,n}, \mu)|] = 0.$

The proof of this lemma is entirely parallel to that of Theorem 6.1, so we omit it. This completes the proof of (7.1).

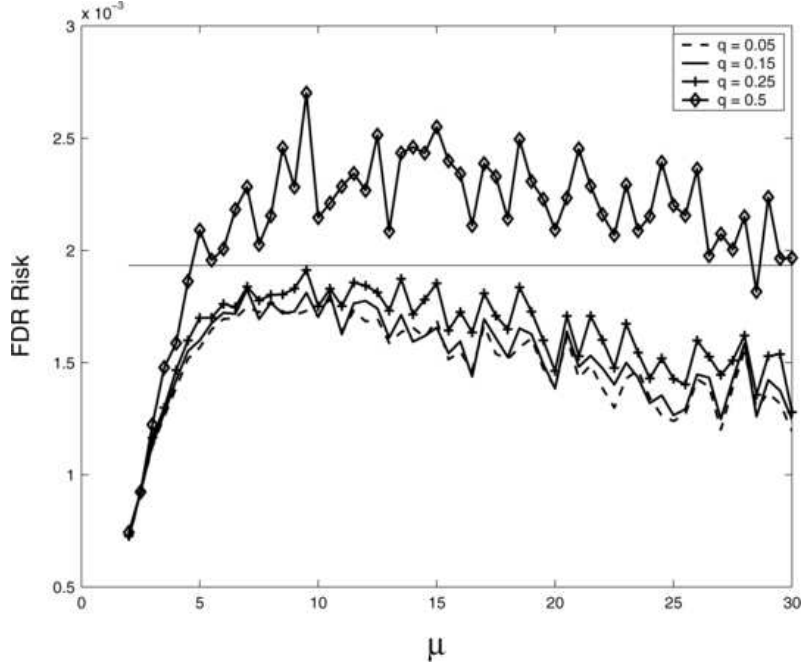

Fig. 4. *Simulation results for FDR thresholding. Curves (dashed, solid, cross, and diamond) describe per-coordinate loss of the FDR procedure with different q values (q = 0.05, 0.15, 0.25) at different two-point mixtures. Here, the mixtures concentrate at 1 and $\mu$ with mass $\varepsilon = \eta/\log(\mu)$ at $\mu$. The horizontal line corresponds to the asymptotic risk expression $\eta \log\log(\frac{1}{\eta})$.*

## 8. Discussion.

8.1. *Illustrations.* We briefly illustrate two key points.

First, we consider finite-sample performance of FDR thresholding. Figure 4 shows the result of FDR thresholding with various values of $q$. It used a sample size $n = 10^6$, sparsity parameters $p = 1$, $\eta = 10^{-3}$ and a range of two-point mixtures of the kind discussed in Theorem 5.1. The figure compares the actual risk of the FDR procedure under a range of situations with the asymptotic limit given by Theorem 1.3. Clearly, the risk depends more strongly on $q$ in finite samples than seems called for by the asymptotic expression in Theorem 1.3. In the simulations, the mixtures were based on various $(\varepsilon, \mu)$ pairs with $\mu$ ranging between 2 and 30 and where, for each $\mu$, $\varepsilon = \frac{\eta}{\log(\mu)}$.

For each $q \in \{0.05, 0.15, 0.25, 0.5\}$, we applied the FDR thresholding estimator $\hat{\mu}_{q,n}^{FDR}$, obtaining an empirical risk measure

$$\hat{R}(q, \mu) = \hat{R}(q, \mu; \eta, n) = \frac{1}{n} \|\log \hat{\mu}_{q,n} - \log \mu\|_2^2.$$

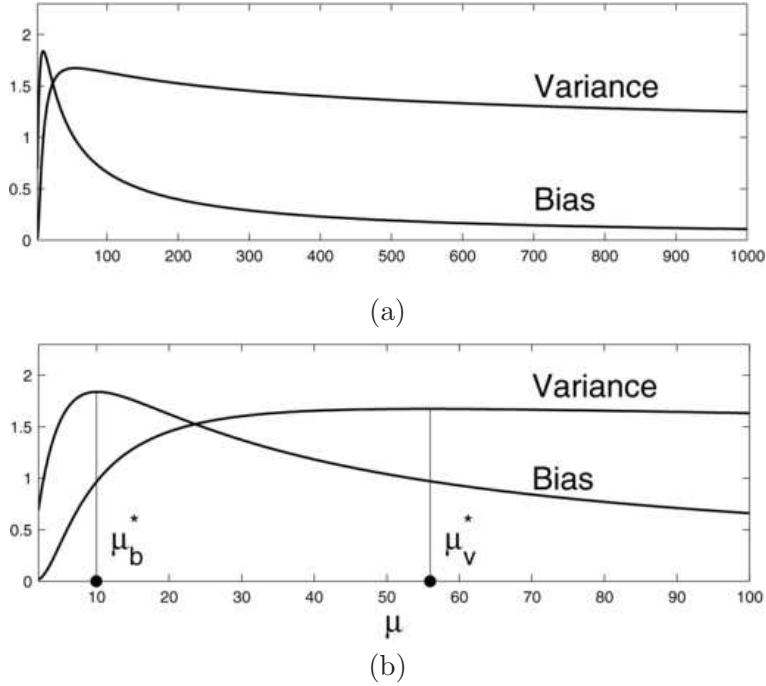

FIG. 5. *Panel (a): The 'bias proxy' $\tilde{B}^2(T_q, G_{\varepsilon,\mu})$ and the 'variance proxy' $\tilde{V}(T_q, G_{\varepsilon,1,\mu})$. Panel (b): Enlargement of (a). The maxima of $\tilde{B}^2(T_q, G_{\varepsilon,\mu})$ and $\tilde{V}(T_q, G_{\varepsilon,\mu})$ are obtained roughly at $\mu_b^*$ and $\mu_v^*$, respectively, with $\mu_b^* = \log(\frac{1}{\eta})/\log\log\log(\frac{1}{\eta})$, $\mu_v^* = \log(\frac{1}{\eta}) \cdot \log\log(\frac{1}{\eta})$. In this figure, $\eta = 10^{-6}$.*

Figure 4 plots $\hat{R}(q, \mu; \eta, n)$ versus $\mu$ for each $q$. As $\mu$ varies between 2 and 30, the empirical FDR risk first increases to a maximum, then decreases; this fits well with our theory. We also note that for $q$ smaller than $1/2$, the empirical FDR risk is not larger than $\eta \log \log(\frac{1}{\eta})$ and when $q$ is close to $1/2$, though the empirical FDR risk can be larger than $\eta \log \log(\frac{1}{\eta})$, it is rarely larger than, say, $1.3 \cdot \eta \log \log(\frac{1}{\eta})$.

Second, we illustrate the behavior of the ideal risk function introduced in the second part of Theorem 5.1. Figure 5 illustrates an example of the ideal risk decomposition into bias proxy and variance proxy, showing the maxima of each and the different ranges over which the two assume their large values.

8.2. *Generalizations.* The approach described here can be directly extended to other settings. Jin has recently derived, by similar methods, asymptotic minimaxity of FDR thresholding for sparse Poisson means obeying $\mu \geq 1$, with most $\mu_i = 1$. This could be useful in situations where we have a collection of 'cells' and expect one event per cell in typical cases, with occasional 'hot spots' containing more than one event per cell.

Preliminary calculations show that a wide range of non-Gaussian additive noises can also be handled by these methods. To see why, note that due to the use of $\log(\mu_i)$ in both the loss measure and parameter set, the results of this paper can be considered a study of FDR thresholding in a situation with *additive* noise having a standard Gumbel distribution. Thus, defining $Y_i = \log(X_i)$, the model of Section 1 posits effectively that

$$Y_i = \theta_i + Z_i, \qquad i = 1, \ldots, n,$$

where, for $\theta_i \geq 0$,

$$\frac{1}{n}\left(\sum_i \theta_i^p\right) \leq \eta^p,$$

we measure loss by $\sum_i (\hat{\theta}_i - \theta_i)^2$ and the noise $Z_i$ obeys $e^{Z_i} \sim \text{Exp}(1)$. Although we have focused on the *one-sided* problem in which $\theta_i \geq 0$ for all $i$, we can certainly generalize the study to treat the *two-sided* problem with $\frac{1}{n}(\sum_i |\theta_i|^p) \leq \eta^p$ and where both $\theta_i > 0$ and $\theta_i < 0$ are possible. Other additive non-Gaussian noises which have been considered include double-exponential. Of course, in considering non-Gaussian distributions, the effectiveness of thresholding depends on the tails of the noise distribution being sufficiently light. Thus, asymptotic minimaxity of thresholding would be doubtful for additive Cauchy noise.

Another generalization concerns dependent settings. In principle, FDR thresholding can still be 'estimating' the FDR functional in large samples, even without i.i.d. stochastic disturbances. Suppose that the $X_i$ are weakly dependent, in such a way that their empirical cdf still converges at a root-$n$ rate. Then all of the above analysis can be carried through in detail without essential change.

One frequently raised question whether the study here could easily be generalized to other distributional settings such as other exponential families. Unfortunately, the results in this paper depend on some properties of the exponential distribution which other exponential families might not have. The most important is the *monotone likelihood ratio* of the family of exponential density functions $\{f_\mu(x), 0 < \mu < \infty : f_\mu(x) = \frac{1}{\mu} e^{-x/\mu} \cdot \mathbb{1}_{\{x>0\}}\}$ [14]; this seems crucial for our argument [12], but some exponential families are not MLR. Jin's study shows that the behavior of the FDR functional in the discrete Poisson setting is essentially different from that of a continuous setting (Gaussian, exponential, etc.). Another frequently raised issue concerns the possibility of working on the original scale instead of the log-scale. However, this does not give rise to a meaningful problem; if we used $\ell^2$-loss on $\mu$ instead of on $\log \mu$, the minimax risk would be infinite.

8.3. *Relation to other work.* There are two points of contact with earlier literature. The first, of course, is with the work of Abramovich, Benjamini, Donoho and Johnstone [3]. Like the present work, [3] proves an asymptotic minimaxity property for the FDR thresholding estimator only for Gaussian noise, and for a subtly different notion of sparsity. In [3], the sparsity parameter $\eta = \eta_n$ so that the sparsity is linked to sample size, which makes sense in a variety of nonparametric estimation applications such as like wavelet denoising [1, 2, 7, 8]. In our work, $\eta$ goes to zero only *after* $n \to \infty$. This simplifies our analysis; the underlying tools in [3]—empirical processes, moderate deviations—are more delicate to deploy than ours. The advantage of our approach seems to lie principally in the ease of generalization to a wider range of non-Gaussian and dependent situations.

The second connection is with the work of Genovese and Wasserman [10]. While they do not consider our multiparameter estimation problem, they do use a Bayesian viewpoint related to Sections 4–6 of our paper. Our approach considers, of course, a different class of Bayesian examples and a different notion of estimation risk. Their paper seems focused on developing intuition and a broader understanding of the FDR approach, while ours uses FDR to attack a specific optimal estimation problem.

**Acknowledgments.** The authors would like to thank Yoav Benjamini and Iain Johnstone for extensive discussions, references and encouragement. These results were reported in J.J.'s PhD Thesis. J.J. thanks committee members for valuable suggestions and encouragement.

DEPARTMENT OF STATISTICS
STANFORD UNIVERSITY
STANFORD, CALIFORNIA 94305
USA

DEPARTMENT OF STATISTICS
PURDUE UNIVERSITY
150 N. UNIVERSITY STREET
WEST LAFAYETT, INDIANA 47907
USA
E-MAIL: jinj@stat.purdue.edu